\documentclass[a4paper,12pt]{article}

\usepackage{latexsym}
\usepackage{amssymb}
\usepackage{theorem}
\usepackage{amsmath}
\usepackage{amscd}
\usepackage{graphicx}
\usepackage{color}

\newtheorem{theorem}{Theorem}[section]
\newtheorem{corollary}[theorem]{Corollary}
\newtheorem{lemma}[theorem]{Lemma}
\newtheorem{example}[theorem]{Example}
\newtheorem{proposition}[theorem]{Proposition}
\newtheorem{remark}[theorem]{Remark}
\newtheorem{definition}[theorem]{Definition}

\newcommand{\demo}{\par\noindent{\it Proof. \/}\ }
\newcommand{\enD}{\hfill $\Box$\vspace{3truemm} \par}
\newcommand{\R}{\mathbb{R}}

\newcommand{\sign}{\operatorname{sign}}
\newcommand{\bn}{\mbox{\boldmath $n$}}
\newcommand{\bt}{\mbox{\boldmath $t$}}
\newcommand{\ba}{\mbox{\boldmath $a$}}
\newcommand{\bb}{\mbox{\boldmath $b$}}

\newcommand{\bu}{\mbox{\boldmath $u$}}
\newcommand{\bv}{\mbox{\boldmath $v$}}
\newcommand{\bw}{\mbox{\boldmath $w$}}

\begin{document}

\title{Bertrand Legendre curves in the unit tangent bundle over Euclidean plane}

\author{Nozomi Nakatsuyama and Masatomo Takahashi}

\date{\today}

\maketitle

\begin{abstract}
We investigate not only the associated curves of regular plane curves, but also those of Legendre curves. 
As associated curves, we consider Bertrand regular plane curves and Bertrand Legendre curves.
These curves contain parallel, evolute and involute curves, as well as evolutoid and involutoid curves. 
Since associated curves may have singular points even if the original curve is regular, Legendre curves provide a suitable framework for investigating the properties of such curves. 
We give existence conditions of Bertrand regular plane curves and Bertrand Legendre curves.
Moreover, we give an inverse operation for Bertrand Legendre curves. 
Furthermore, we define a mapping between sets of Legendre curves using Bertrand Legendre curves and prove that this mapping is bijective up to equivalence relations. 
\end{abstract}

\renewcommand{\thefootnote}{\fnsymbol{footnote}}

%%%%%%%%%%%%%%%% Section 1 %%%%%%%%%%%%%%%%%%
\section{Introduction}

The Bertrand and Mannheim curves are classical objects in differential geometry (\cite{Aminov, Banchoff-Lovett, Berger-Gostiaux, Bertrand, doCarmo, Honda-Takahashi-2020, HCIP, Izumiya-Takeuchi1, Kuhnel, Liu-Wang, Papaioannou-Kiritsis, Struik}). 
The Bertrand (respectively, Mannheim) curve is a space curve whose principal normal line is the same as the principal normal (respectively, bi-normal) line of another curve. 
We investigated Bertrand type curves as general cases in Euclidean 3-space (cf. \cite{Nakatsuyama, Nakatsuyama-Takahashi2}). 
In this paper, we apply the idea to the plane curves. 
We investigate not only the associated curves of regular plane curves, but also those of Legendre curves. 
As associated curves, we consider Bertrand regular plane curves and Bertrand Legendre curves. 
These curves contain parallel, evolute and involute curves, as well as evolutoid and involutoid curves. 
Since associated curves may have singular points even if the original curve is regular, Legendre curves provide a suitable framework for investigating the properties of such curves. 
A Legendre curve in the unit tangent bundle over Euclidean plane is a plane curve with a moving frame. 
Then we can define the (Legendre) curvature of the Legendre curve  (cf. \cite{Fukunaga-Takahashi-2013}). 
The curvature is a complete invariant for Legendre curves up to congruence (Euclidean motion) like as the curvature of regular plane curves. 
The existence and uniqueness theorems for the curvature are valid (cf. \cite{Fukunaga-Takahashi-2013}). 
The definitions of the evolutoid of regular plane curves and frontals have already investigated in \cite{Giblin-Warder, Izumiya-Takeuchi-2019}, and the definition of the involutoid (tanvolute) of regular plane curves has already investigated in \cite{AM}. 
Moreover, the evolutoid and involutoid of spherical Legendre curves investigated in \cite{Li-Pei}.
However, as far as we know, the explicit form of the definition of the involutoid of regular plane curves and of Legendre curves (frontals)  in the unit tangent bundle over Euclidean plane cannot be found. 
\par
We give the conditions for the existence of Bertrand regular plane curves and Bertrand Legendre curves (Theorems \ref{(v,w)-Bertrand-regular} and \ref{vw-Legendre}). 
Then we give the explicit form of the involutoid. 
We also give new correspondences which connects between evolutes and involutes of Legendre curves.  
Moreover, we give an inverse operation for Bertrand Legendre curves (Theorem \ref{inverse-relation}). 
Furthermore, we define a mapping between sets of Legendre curves using Bertrand Legendre curves and prove that this mapping is bijective up to equivalence relations  (Theorem \ref{bijective}).
\par
We shall assume throughout the whole paper that all maps and manifolds are $C^{\infty}$ unless the contrary is explicitly stated.

%%%%%%%%%%%%%%%%%%%%%%%%%%%%%%%%%%%%%%%%%%%%
\bigskip
\noindent
{\bf Acknowledgements}. 
The first author was supported by JST SPRING Grant Number JPMJSP2153. 
The second author was partially supported by JSPS KAKENHI Grant Number JP 24K06728.

%%%%%%%%%%%%%%%%% Section 2 %%%%%%%%%%%%%%%%%
\section{Preliminaries}

Let $\R^2$ be the $2$-dimensional Euclidean space equipped with the inner product $\ba \cdot \bb = a_1 b_1 + a_2 b_2$, 
where $\ba = (a_1, a_2)$ and $\bb = (b_1, b_2) \in \R^2$. 
The norm of $\ba$ is given by $\vert \ba \vert = \sqrt{\ba \cdot \ba}$.
We denote the anti-clockwise rotation of $\pi/2$ on $\R^2$ by $J$.
Let $S^1$ be the unit circle in $\R^2$, that is, $S^1=\{\ba \in \R^2| |\ba|=1\}$.
We review the theories of regular curves (cf. \cite{Banchoff-Lovett, doCarmo, Gibson, Gray}) and Legendre curves (cf. \cite{Fukunaga-Takahashi-2013, Fukunaga-Takahashi-2015}).
%%%%%
\subsection{Regular plane curves}

Let $I$ be an interval of $\R$ and let $\gamma:I \to \R^2$ be a regular plane curve, that is, $\dot{\gamma}(t) \not=0$ for all $t \in I$, where $\dot{\gamma}(t)=(d\gamma/dt)(t)$. 
If we take the arc-length parameter $s$, that is, $|\gamma'(s)|=1$ for all $s$, then the tangent vector, the normal vector are given by $\bt(s)=\gamma'(s), \ \bn(s)=J(\bt(s)),$ where $\gamma'(s)=(d\gamma/ds)(s)$. 
Then $\{\bt(s),\bn(s)\}$ is a moving frame of $\gamma(s)$ and we have the Frenet formula: 
$$
\left(
\begin{array}{c}
\bt'(s)\\
\bn'(s)
\end{array}
\right)
=
\left(
\begin{array}{cc}
0&\kappa(s)\\
-\kappa(s)&0
\end{array}
\right)
\left(
\begin{array}{c}
\bt(s)\\
\bn(s)
\end{array}
\right),
$$
where $\kappa(s)=\bt'(s) \cdot \bn(s).$
If we take a general parameter $t$, then the tangent vector, the normal vector are given by
$\bt(t)={\dot{\gamma}(t)}/{|\dot{\gamma}(t)|}, \bn(t)=J(\bt(t)).$
Then $\{\bt(t),\bn(t)\}$ is a moving frame of $\gamma(t)$ and we have the Frenet formula: 
$$
\left(
\begin{array}{c}
\dot{\bt}(t)\\
\dot{\bn}(t)
\end{array}
\right)
=
\left(
\begin{array}{cc}
0&|\dot{\gamma}(t)|\kappa(t)\\
-|\dot{\gamma}(t)|\kappa(t)&0
\end{array}
\right)
\left(
\begin{array}{c}
\bt(t)\\
\bn(t)
\end{array}
\right),
$$
where 
$\kappa(t)={{\rm det}(\dot{\gamma}(t),\ddot{\gamma}(t))}/{|\dot{\gamma}(t)|^3}.$

%%%%%
\subsection{Legendre curves}

Let $(\gamma,\nu):I \to \R^2 \times S^1$ be a smooth mapping.

%%%%%
\begin{definition}{\rm
We say that $(\gamma,\nu):I \to \R^2 \times S^1$ is a {\it Legendre curve} if $\dot{\gamma}(t) \cdot \nu(t)=0$ for all $t \in I$. 
We also say that $\gamma$ is a {\it frontal} if there exists a smooth map $\nu: I \to S^1$ such that $(\gamma,\nu)$ is a Legendre curve. 
}
\end{definition}
%%%%% 
We denote $\mu:I \to S^1$ by $\mu(t)=J(\nu(t))$.
Then $\{\nu(t),\mu(t) \}$ is a moving frame of $\gamma(t)$. 
Therefore, we have the following Frenet formula:
\begin{align*}
\left(
\begin{array}{c}
\dot{\nu}(t)\\
\dot{\mu}(t)
\end{array}
\right)
=
\left(
\begin{array}{cc}
0&\ell(t)\\
-\ell(t)&0
\end{array}
\right)
\left(
\begin{array}{c}
\nu(t)\\
\mu(t)
\end{array}
\right), \ 
\dot{\gamma}(t)
=\beta(t)\mu(t),
\end{align*}
where $\ell(t)=\dot{\nu}(t) \cdot \mu(t), \beta(t)=\dot{\gamma}(t) \cdot \mu(t)$.
The mapping $(\ell,\beta):I \to \R^2$ is called the (Legendre) curvature of $(\gamma,\nu)$. 
By definition, $t \in I$ is a singular point of $\gamma$, that is, $\dot{\gamma}(t)=0$ if and only if $\beta(t)=0$.
We say that $t \in I$ is an {\it inflection point} of $\gamma$ (or, $(\gamma,\nu)$) if $\ell(t)=0$. 
\par
Let $(\gamma,\nu):I \to \R^2 \times S^1$ be a Legendre curve with curvature $(\ell,\beta)$. 
Then $(\gamma,-\nu):I \to \R^2 \times S^1$ is also a Legendre curve with curvature $(\ell,-\beta)$.

Obviously, a regular curve is a frontal. 
We give a relation between curvatures of Legendre curves and of regular curves. 
%%%%%
\begin{proposition}\label{regular-Legendre}
$(1)$ Let $(\gamma,\nu):I \to \R^2 \times S^1$ be a Legendre curve with $(\ell,\beta)$.
If $\gamma$ is a regular plane curve with curvature $\kappa$, then we have $\ell=\kappa |\beta|$.
\par
$(2)$
If $\gamma:I \to \R^2$ is a regular plane curve with curvature $\kappa$, then $(\gamma,\bn):I \to \R^2 \times S^1$ is a Legendre curve with curvature $(\ell,\beta)=(|\dot{\gamma}|\kappa, -|\dot{\gamma}|)$.
\end{proposition}
%%%%%

%$\bullet$ Definition of type of singular points.

We say that a smooth plane curve $\gamma:(I,t_0) \to \R^2$ at $t_0$ is a $j/i$-cusp (or, $(i,j)$-cusp), where $(i, j)= (2, 3),(2, 5)$, $(3, 4)$ and $(3, 5)$ if $\gamma$ is $ \mathcal{A}$-equivalent to the germ $t \mapsto(t^i,t^j)$ at the origin. 
For curves with $j/i$-cusps in $\R^2$, the criteria are known in \cite{Bruce-Gaffney,Porteous}.

%$\bullet$ Criteria of singular points (3/2, 4/3, 5/2, 5/3 cusp) as Proposition.
%%%%%
\begin{proposition}\label{singular-point:criteria}
Let $(\gamma,\nu):I \to \R^2 \times S^1$ be a Legendre curve with $(\ell,\beta)$.
\par
$(1)$ $\gamma$ has a $3/2$-cusp at $t_0$ if and only if $\beta(t_0)=0$, $\dot\beta(t_0)\neq0$ and $\ell(t_0)\neq0$.
\par
$(2)$ $\gamma$ has a $5/2$-cusp at $t_0$ if and only if $\beta(t_0)=\ell(t_0)=0$, $\dot\beta(t_0)\neq0$ and $\ddot\ell(t_0)\dot\beta(t_0)-\dot\ell(t_0)\ddot\beta(t_0)\neq0$.
\par
$(3)$ $\gamma$ has a $4/3$-cusp at $t_0$ if and only if $\beta(t_0)=\dot\beta(t_0)=0$, $\ddot\beta(t_0)\neq0$ and $\ell(t_0)\neq0$.
\par
$(4)$ $\gamma$ has a $5/3$-cusp at $t_0$ if and only if $\beta(t_0)=\dot\beta(t_0)=\ell(t_0)=0$, $\ddot\beta(t_0)\neq0$ and $\dot\ell(t_0)\neq0$.
\end{proposition}
%%%%%
%%%%%
\section{Bertrand regular plane curves}

Let $\gamma, \overline{\gamma}: I \to \R^2$ be regular plane curves and $\bv, \overline{\bw}: I \to S^1$ be smooth curves. 

%%%%%
\begin{definition}\label{Bertrand-regular}
{\rm 
We say that $\gamma$ and $\overline{\gamma}$ are {\it $(\bv,\overline{\bw})$-regular mates} if there exists a smooth function $\lambda:I \to \R$ with $\lambda \not\equiv 0$ such that $\overline{\gamma}(t)=\gamma(t)+\lambda(t)\bv(t)$ and $\bv(t)=\overline{\bw}(t)$ for all $t \in I$. 
We also say that $\gamma:I \to \R^2$ is a {\it $(\bv,\overline{\bw})$-Bertrand regular (plane) curve} if there exists a regular plane curve $\overline{\gamma}:I \to \R^2$ such that $\gamma$ and $\overline{\gamma}$ are $(\bv,\overline{\bw})$-regular mates. 
}
\end{definition}
%%%%%

We clarify the notation $\lambda \not\equiv 0$. 
Throughout this paper, $\lambda \not\equiv 0$ means that $\{t \in I | \lambda(t) \not=0\}$ is a dense subset of $I$. 
Then $\lambda$ is not identically zero for any non-trivial subintervals of $I$. 
It follows that $\gamma$ and $\overline{\gamma}$ are different plane curves for any non-trivial subintervals of $I$. 
Note that if $\lambda$ is constant, then $\lambda \not\equiv 0$ means that $\lambda$ is a non-zero constant. 
\par
We give a characterization of Bertrand regular curves. 
Let $\gamma: I \to \R^2$ be a regular plane curve with curvature $\kappa$.
By a parameter change, we may assume that $s$ is the arc-length parameter of $\gamma$. 
Note that $s$ is not the arc-length parameter of $\overline{\gamma}$.
Since $\bv, \overline{\bw}: I \to S^1$ are smooth curves, there exist smooth functions $\theta,\tau:I \to \R$ such that $\bv(s)=\cos \theta(s) \bt(s)+\sin \theta(s) \bn(s)$ and $\overline{\bw}(s)=\cos \tau(s) \overline{\bt}(s)+\sin \tau(s) \overline{\bn}(s)$. 

%%%%%
\begin{theorem}\label{(v,w)-Bertrand-regular}
Let $\bv,\overline{\bw}:I \to S^1$ be smooth curves and $\gamma:I \to \R^2$ be a regular plane curve with curvature $\kappa$. 
Then $\gamma:I \to \R^2$ is a $(\bv,\overline{\bw})$-Bertrand regular curve if and only if there exist smooth functions $\lambda, \theta, \tau:I \to \R$ with $\lambda \not\equiv 0$ such that 
\begin{align}
&(\cos \theta(s)+{\lambda}'(s))\sin \tau(s)-(\sin \theta(s)-\lambda(s)({\theta}'(s)+\kappa(s)))\cos \tau(s)=0, \label{regular-mate-condition} \\
&(\cos \theta(s)+{\lambda}'(s))\cos \tau(s)+(\sin \theta(s)-\lambda(s)({\theta}'(s)+\kappa(s)))\sin \tau(s) \not=0 \label{regular-condition}
\end{align} 
for all $s \in I$.
\end{theorem}
%%%%%
\demo 
Suppose that $\gamma:I \to \R^2$ is a $(\bv,\overline{\bw})$-Bertrand regular curve.
Then there exists a smooth function $\lambda:I \to \R$ with $\lambda \not\equiv 0$ such that $\overline{\gamma}(s)=\gamma(s)+\lambda(s)\bv(s)$ and $\bv(s)=\overline{\bw}(s)$ for all $s \in I$, where there exist smooth functions $\theta, \tau:I \to \R$ such that $\bv(s)=\cos \theta(s)\bt(s)+\sin \theta(s)\bn(s)$ and $\overline{\bw}(s)=\cos \tau(s) \overline{\bt}(s)+\sin \tau(s) \overline{\bn}(s)$. 
By differentiating $\overline{\gamma}(s)=\gamma(s)+\lambda(s)\bv(s)$, we have 
\begin{align*}
|\dot{\overline{\gamma}}(s)|\overline{\bt}(s)=&\bt(s)+\lambda'(s)(\cos\theta(s)\bt(s)+\sin\theta(s)\bn(s))\\
&+\lambda(s)(-(\theta'(s)+\kappa(s))\sin\theta(s)\bt(s)+(\theta'(s)+\kappa(s))\cos\theta(s)\bn(s)).
%=&(1+\lambda'(s)\cos\theta(s)-\lambda(s)(\theta'(s)+\kappa(s))\sin\theta(s))\bt(s)\\
%&+(\lambda'(s)\sin\theta(s)+\lambda(s)(\theta'(s)+\kappa(s))\cos\theta(s))\bn(s).
\end{align*}
By the above notations, we have 
$$
\begin{pmatrix}
J({\bv}(s)) \\
{\bv}(s)
\end{pmatrix}
=
\begin{pmatrix}
\cos \theta(s) & -\sin \theta(s)\\
\sin \theta(s) & \cos \theta(s)
\end{pmatrix}
\begin{pmatrix}
{\bn}(s) \\
{\bt}(s)
\end{pmatrix}, \
\begin{pmatrix}
J(\overline{\bw}(s)) \\
\overline{\bw}(s)
\end{pmatrix}
=
\begin{pmatrix}
\cos \tau(s) & -\sin \tau(s)\\
\sin \tau(s) & \cos \tau(s)
\end{pmatrix}
\begin{pmatrix}
\overline{\bn}(s) \\
\overline{\bt}(s)
\end{pmatrix}.
$$
Moreover, since $\bv(s)=\overline{\bw}(s)$, we have $J(\bv(s))=J(\overline{\bw}(s))$, that is, 
$$
\cos\theta(s)\bn(s)-\sin\theta(s)\bt(s)=\cos\tau(s)\overline{\bn}(s)-\sin\tau(s)\overline{\bt}(s).
$$ 
Therefore, we have
\begin{align*}
|\dot{\overline\gamma}(s)|\overline{\bt}(s)&=(\cos\theta(s)+\lambda'(s))\bv(s)-(\sin\theta(s)-\lambda(s)(\theta'(s)+\kappa(s)))J(\bv(s))\\
%&=(\cos\theta(s)+\lambda'(s))(\cos \theta(s) \bt(s)+\sin \theta(s) \bn(s))\\
%&\quad -(\sin\theta(s)-\lambda(s)(\theta'(s)+\kappa(s)))(\cos\theta(s)\bn(s)-\sin\tau(s)\bt(s))\\
&=(\cos\theta(s)+\lambda'(s))(\cos \tau(s) \overline{\bt}(s)+\sin \tau(s) \overline{\bn}(s))\\
&\quad -(\sin\theta(s)-\lambda(s)(\theta'(s)+\kappa(s)))(\cos\tau(s)\overline{\bn}(s)-\sin\tau(s)\overline{\bt}(s))\\
&=\left((\cos\theta(s)+\lambda'(s))\cos \tau(s)+(\sin\theta(s)-\lambda(s)(\theta'(s)+\kappa(s)))\sin\tau(s) \right)\overline{\bt}(s)\\
&\quad +\left((\cos\theta(s)+\lambda'(s))\sin \tau(s)-(\sin\theta(s)-\lambda(s)(\theta'(s)+\kappa(s)))\cos\tau(s) \right)\overline{\bn}(s).
\end{align*}
It follows that we have equations (\ref{regular-mate-condition}) and (\ref{regular-condition}).
\par
Conversely, suppose that equations (\ref{regular-mate-condition}) and (\ref{regular-condition}) satisfy. 
Let $\overline{\gamma}:I \to \R^2$ be $\overline{\gamma}(t)=\gamma(t)+\lambda(t)\bv(t)$, where $\bv(s)=\cos \theta(s) \bt(s)+\sin \theta(s) \bn(s)$. 
By a direct calculation, we have 
\begin{align*}
\dot{\overline\gamma}(s)=&\left((\cos\theta(s)+\lambda'(s))\cos \tau(s)+(\sin\theta(s)-\lambda(s)(\theta'(s)+\kappa(s)))\sin\tau(s) \right)\overline{\bt}(s)\neq0.
\end{align*}
It follows that $\overline{\gamma}$ is a regular plane curve. 
Moreover, since 
\begin{align*}
|\dot{\overline\gamma}(s)|\overline{\bt}(s) &=(1+\lambda'(s)\cos\theta(s)-\lambda(s)(\theta'(s)+\kappa(s))\sin\theta(s))\bt(s) \\
&\quad +(\lambda'(s)\sin\theta(s)+\lambda(s)(\theta'(s)+\kappa(s))\cos\theta(s))\bn(s),
\end{align*}
then we have 
\begin{align*}
%|\dot{\overline\gamma}(s)|J(\overline{\bt}(s))=& \ (1+\lambda'(s)\cos\theta(s)-\lambda(s)(\theta'(s)+\kappa(s))\sin\theta(s))J(\bt(s))\\
%&+(\lambda'(s)\sin\theta(s)+\lambda(s)(\theta'(s)+\kappa(s))\cos\theta(s))J(\bn(s))\\
|\dot{\overline\gamma}(s)|\overline{\bn}(s)=& \ (1+\lambda'(s)\cos\theta(s)-\lambda(s)(\theta'(s)+\kappa(s))\sin\theta(s))\bn(s)\\
&-(\lambda'(s)\sin\theta(s)+\lambda(s)(\theta'(s)+\kappa(s))\cos\theta(s))\bt(s).
\end{align*}
Therefore, we have 
\begin{align*}
\overline{\bt}(s)=& \ \frac{1+\lambda'(s)\cos\theta(s)-\lambda(s)(\theta'(s)+\kappa(s))\sin\theta(s)}{(\cos\theta(s)+\lambda'(s))\cos \tau(s)+(\sin\theta(s)-\lambda(s)(\theta'(s)+\kappa(s)))\sin\tau(s)}\bt(s)\\
&+\frac{\lambda'(s)\sin\theta(s)+\lambda(s)(\theta'(s)+\kappa(s))\cos\theta(s)}{(\cos\theta(s)+\lambda'(s))\cos \tau(s)+(\sin\theta(s)-\lambda(s)(\theta'(s)+\kappa(s)))\sin\tau(s)}\bn(s),\\
\overline{\bn}(s)=& \ -\frac{\lambda'(s)\sin\theta(s)+\lambda(s)(\theta'(s)+\kappa(s))\cos\theta(s)}{(\cos\theta(s)+\lambda'(s))\cos \tau(s)+(\sin\theta(s)-\lambda(s)(\theta'(s)+\kappa(s)))\sin\tau(s)}\bt(s)\\
&+\frac{1+\lambda'(s)\cos\theta(s)-\lambda(s)(\theta'(s)+\kappa(s))\sin\theta(s)}{(\cos\theta(s)+\lambda'(s))\cos \tau(s)+(\sin\theta(s)-\lambda(s)(\theta'(s)+\kappa(s)))\sin\tau(s)}\bn(s).\\
\end{align*}
By a direct calculation, we have 
$$
\overline{\bw}(s)=\cos\tau(s)\overline{\bt}(s)+\sin\tau(s)\overline{\bn}(s)=\cos\theta(s)\bt(s)+\sin\theta(s)\bn(s)=\bv(s).
$$ 
It follows that $\gamma:I \to \R^2$ is a $(\bv,\overline{\bw})$-Bertrand regular curve.
\enD
%%%%%
%%%%%
\begin{remark}{\rm 
In general, a deformation of $\gamma$ may have singular points.
The condition \eqref{regular-condition} is nothing but a regular condition of $\overline{\gamma}$. 
The curvature of $\overline{\gamma}$ is given by 
$$
\overline{\kappa}(s)=\frac{\theta'(s)-\tau'(s)+\kappa(s)}{|(\cos\theta(s)+\lambda'(s))\cos \tau(s)+(\sin\theta(s)-\lambda(s)(\theta'(s)+\kappa(s)))\sin\tau(s)|}.
$$
Hence, $\cos \theta(s)+\lambda'(s)=0$ and $\sin \theta(s)-\lambda(s)(\theta'(s)+\kappa(s))=0$ if and only if $s$ is a singular point of $\overline{\gamma}$. 
}
\end{remark}
%%%%%
\par
By Theorem \ref{(v,w)-Bertrand-regular}, we have the special cases of Bertrand regular curves.
%%%%%
\begin{corollary}\label{special-cases}
Under the above notations, we have the following.
\par
$(1)$ Suppose that $\theta(s)=0$ and $\tau(s)=0$ for all $s \in I$. 
Then $\gamma$ is a $(\bt,\overline{\bt})$-Bertrand regular curve if and only if $\kappa(s)=0$ and $1+\lambda'(s) \not=0$ for all $s \in I$. 
It follows that $\overline{\gamma}$ and $\gamma$ are a part of a line.
\par
$(2)$ Suppose that $\theta(s)=0$ and $\tau(s)=\pi/2$ for all $s \in I$. 
Then $\gamma$ is a $(\bt,\overline{\bn})$-Bertrand regular curve if and only if $\lambda(s)=-s+c$, where $c$ is a constant and $\kappa(s) \not=0$ for all $s \in I$. 
It follows that $\overline{\gamma}$ is an involute of $\gamma$.
\par
$(3)$ Suppose that $\theta(s)=\pi/2$ and $\tau(s)=0$ for all $s \in I$. 
Then $\gamma$ is an $(\bn,\overline{\bt})$-Bertrand regular curve if and only if $1-\kappa(s)\lambda(s)=0$ and $\kappa'(s) \not=0$ for all $s \in I$. 
It follows that $\overline{\gamma}$ is an evolute of $\gamma$.
\par
$(4)$ Suppose that $\theta(s)=\pi/2$ and $\tau(s)=\pi/2$ for all $s \in I$. 
Then $\gamma$ is an $(\bn,\overline{\bn})$-Bertrand regular curve if and only if $\lambda(s)=\lambda$ is a constant and $1-\lambda \kappa(s) \not=0$ for all $s \in I$. 
It follows that $\overline{\gamma}$ is a parallel curve of $\gamma$.
\par
$(5)$ Suppose that $\theta(s)$ is a constant $\theta$ and $\tau(s)=0$ for all $s \in I$. 
Then $\gamma$ is a $(\bv,\overline{\bt})$-Bertrand regular curve if and only if $\sin \theta-\lambda(s)\kappa(s)=0$ and $\cos \theta+\lambda'(s) \not=0$ for all $s \in I$. 
It follows that $\overline{\gamma}$ is an evolutoid of $\gamma$.
\par
$(6)$ Suppose that $\theta(s)=0$ and $\tau(s)$ is a constant $\tau$ for all $s \in I$. 
Then $\gamma$ is a $(\bt,\overline{\bw})$-Bertrand regular curve if and only if $(1+\lambda'(s)) \sin \tau +\lambda(s)\kappa(s)\cos \tau=0$ and $(1+\lambda'(s)) \cos \tau -\lambda(s)\kappa(s)\sin  \tau \not=0$ for all $s \in I$. 
It follows that $\overline{\gamma}$ is an involutoid of $\gamma$.
\end{corollary}
%%%%%

%%%%%%%%%%%%%%%%% Section 4 %%%%%%%%%%%%%%%%%
\section{Bertrand Legendre curves}

Let $(\gamma,\nu), (\overline{\gamma},\overline{\nu}):I \to \R^2 \times S^1$ be Legendre curves and $\bv,\overline{\bw}: I \to S^1$ be smooth curves. 
Since $\bv,\overline{\bw}: I \to S^1$ are smooth curves, there exist smooth functions $\theta,\tau:I \to \R$ such that $\bv(t)=\cos \theta(t) \nu(t)+\sin \theta(t) \mu(t)$ and $\overline{\bw}(t)=\cos \tau(t) \overline{\nu}(t)+\sin \tau(t) \overline{\mu}(t)$. 
Then we say that $\bv,\overline{\bw}:I \to S^1$ are mappings with $\theta,\tau: I \to \R$. 

%%%%%
\begin{definition}\label{Bertrand-type-Legendre}{\rm
We say that $(\gamma,\nu)$ and $(\overline{\gamma},\overline{\nu})$ are {\it $(\bv,\overline{\bw})$-mates} if there exists a smooth function $\lambda:I \to \R$ with $\lambda \not\equiv 0$ such that $\overline{\gamma}(t)=\gamma(t)+\lambda(t)\bv(t)$ and $\bv(t)= \overline{\bw}(t)$ for all $t \in I$. 
Then we say that $(\gamma,\nu)$ and $(\overline{\gamma},\overline{\nu})$ are {\it $(\bv,\overline{\bw})$-mates} with $\lambda$. 
We also say that  $(\gamma,\nu)$ is a {\it $(\bv,\overline{\bw})$-Bertrand Legendre curve} if there exists a Legendre curve $(\overline{\gamma},\overline{\nu})$ such that $(\gamma,\nu)$ and $(\overline{\gamma},\overline{\nu})$ are $(\bv,\overline{\bw})$-mates.
}
\end{definition}
%%%%%
\par
%Let $(\gamma,\nu), (\overline{\gamma},\overline{\nu}):I \to \R^2 \times S^1$ be Legendre curves with curvature $(\ell,\beta), (\overline{\ell},\overline{\beta})$, respectively.  
%Since $\bv,\overline{\bw}: I \to S^1$ are smooth curves, there exist $\theta,\tau:I \to \R$ such that $\bv(t)=\cos \theta(t) \nu(t)+\sin \theta(t) \mu(t)$ and $\overline{\bw}(t)=\cos \tau(t) \overline{\nu}(t)+\sin \tau(t) \overline{\mu}(t)$.
We give a characterization of the Bertrand Legendre curve.

%%%%%
\begin{theorem}\label{vw-Legendre}
Let $\bv, \overline{\bw} :I\to S^1$ be smooth curves and $(\gamma,\nu):I \to \R^2 \times S^1$ be a Legendre curve with curvature $(\ell,\beta)$.
Then $(\gamma,\nu):I \to \R^2 \times S^1$ is a $(\bv,\overline{\bw})$-Bertrand Legendre curve if and only if there exist smooth functions $\lambda, \theta, \tau:I \to \R$ with $\lambda \not\equiv 0$ such that
\begin{align}\label{vw-Legendre-condition}
(\beta(t) \sin \theta(t)+\dot{\lambda}(t)) \cos \tau(t)-(\beta(t)\cos \theta(t)+\lambda(t)(\dot{\theta}(t)+\ell(t)))\sin \tau(t)=0
\end{align} 
for all $t \in I$.
\end{theorem}
%%%%%
\demo
Suppose that $(\gamma,\nu):I \to \R^2 \times S^1$ is a $(\bv,\overline{\bw})$-Bertrand Legendre curve.
Then there exist a Legendre curve $(\overline{\gamma},\overline{\nu}):I \to \R^2 \times S^1$ and a smooth function $\lambda:I \to \R$ with $\lambda \not\equiv 0$ such that $\overline{\gamma}(t)=\gamma(t)+\lambda(t)\bv(t)$ and $\bv(t)=\overline{\bw}(t)$ for all $t \in I$, where there exist smooth functions $\theta, \tau:I \to \R$ such that $\bv(t)=\cos \theta(t) \nu(t)+\sin \theta(t) \mu(t)$ and $\overline{\bw}(t)=\cos \tau(t) \overline{\nu}(t)+\sin \tau(t) \overline{\mu}(t)$.  
By differentiating $\overline{\gamma}(t)=\gamma(t)+\lambda(t)\bv(t)$, we have 
\begin{align*}
{\overline{\beta}}(t)\overline{\mu}(t)=&\beta(t)\mu(t)+\dot\lambda(t)\bv(t)+\lambda(t)(\dot\theta(t)-\ell(t))(\cos\theta(t)\mu(t)-\sin\theta(t)\nu(t)).
\end{align*}
By the above notations, we have 
$$
\begin{pmatrix}
J({\bv}(t)) \\
{\bv}(t)
\end{pmatrix}
=
\begin{pmatrix}
\cos \theta(t) & -\sin \theta(t)\\
\sin \theta(t) & \cos \theta(t)
\end{pmatrix}
\begin{pmatrix}
{\mu}(t) \\
{\nu}(t)
\end{pmatrix}, \
\begin{pmatrix}
J(\overline{\bw}(t)) \\
\overline{\bw}(t)
\end{pmatrix}
=
\begin{pmatrix}
\cos \tau(t) & -\sin \tau(t)\\
\sin \tau(t) & \cos \tau(t)
\end{pmatrix}
\begin{pmatrix}
\overline{\mu}(t) \\
\overline{\nu}(t)
\end{pmatrix}.
$$
Moreover, since $\bv(t)=\overline{\bw}(t)$, we have $J(\bv(t))=J(\overline{\bw}(t))$, that is, 
$$
\cos \theta(t)\mu(t)-\sin\theta(t)\nu(t)=\cos\tau(t)\overline{\mu}(t)-\sin\tau(t)\overline{\nu}(t).
$$ 
Therefore, we have
\begin{align*}
\overline{\beta}(t)\overline{\mu}(t)&=(\beta(t)\cos\theta(t)+\lambda(t)(\dot\theta(t)+\ell(t)))J(\bv(t))+(\beta(t)\sin\theta(t)+\dot\lambda(t))\bv(t)\\
&=(\beta(t)\cos\theta(t)+\lambda(t)(\dot\theta(t)+\ell(t)))(\cos\theta(t)\mu(t)-\sin\theta(t)\nu(t))\\
&\quad +(\beta(t)\sin\theta(t)+\dot\lambda(t))(\sin \theta(t) \mu(t)+\cos \theta(t) \nu(t))\\
&=(\beta(t)\cos\theta(t)+\lambda(t)(\dot\theta(t)+\ell(t)))(\cos\tau(t)\overline{\mu}(t)-\sin\tau(t)\overline{\nu}(t))\\
&\quad +(\beta(t)\sin\theta(t)+\dot\lambda(t))(\sin \tau(t) \overline{\mu}(t)+\cos \tau(t) \overline{\nu}(t))\\
&=\left((\beta(t)\cos\theta(t)+\lambda(t)(\dot\theta(t)+\ell(t)))\cos\tau(t)+(\beta(t)\sin\theta(t)+\dot\lambda(t))\sin \tau(t)\right)\overline{\mu}(t)\\
&\quad +\left((\beta(t)\cos\theta(t)+\lambda(t)(\dot\theta(t)+\ell(t)))\sin\tau(t)-(\beta(t)\sin\theta(t)+\dot\lambda(t))\cos \tau(t)\right)\overline{\nu}(t)
\end{align*}
It follows that we have equation (\ref{vw-Legendre-condition}).
\par
Conversely, suppose that equation (\ref{vw-Legendre-condition}) satisfies. 
Let $(\overline{\gamma},\overline{\nu}):I\to\R^2\times S^1$ be $\overline{\gamma}(t)=\gamma(t)+\lambda(t)\bv(t)$ and $\overline{\nu}(t)=\cos(\theta(t)-\tau(t))\nu(t)+\sin(\theta(t)-\tau(t))\mu(t)$, where $\bv(t)=\cos \theta(t) \nu(t)+\sin \theta(t) \mu(t)$. 
Since 
\begin{align*}
\dot{\overline\gamma}(t)&=\beta(t)\mu(t)+\dot\lambda(t)\bv(t)+\lambda(t)(\dot\theta(t)+\ell(t))(\cos\theta(t)\mu(t)-\sin\theta(t)\nu(t))\\
&=(\dot\lambda(t)\cos\theta(t)-\lambda(t)(\dot\theta(t)+\ell(t))\sin\theta(t))\nu(t)\\
&\quad+(\beta(t)+\dot\lambda(t)\sin\theta(t)+\lambda(t)(\dot\theta(t)+\ell(t))\cos\theta(t))\mu(t),
\end{align*}
$\dot{\overline{\gamma}}(t)\cdot\overline{\nu}(t)=0$ for all $t \in I$.
It follows that $(\overline{\gamma},\overline{\nu})$ is a Legendre curve. 
Then $\overline{\mu}(t)=J(\overline{\nu}(t))=-\sin(\theta(t)-\tau(t))\nu(t)+\cos(\theta(t)-\tau(t))\mu(t)$.
%\begin{align*}
%\overline{\beta}(t)\overline{\mu}(t)=&\ \overline{\beta}(t)\left(\cos\tau(t) J(\overline{\bw}(t))+\sin\tau(t)\overline{\bw}(t)\right)\\
%=& \ \left(-(\beta(t)\cos\theta(t)+\lambda(t)(\dot\theta(t)+\ell(t)))\sin\theta(t)+(\beta(t)\sin\theta(t)+\dot\lambda(t))\cos\theta(t)\right)\nu(t)\\
%&+\left((\beta(t)\cos\theta(t)+\lambda(t)(\dot\theta(t)+\ell(t)))\cos\theta(t)+(\beta(t)\sin\theta(t)+\dot\lambda(t))\sin\theta(t)\right)\mu(t),
%\end{align*}
%\begin{align*}
%\overline{\beta}(t)J(\overline{\mu}(t))=& \ \left(-(\beta(t)\cos\theta(t)+\lambda(t)(\dot\theta(t)+\ell(t)))\sin\theta(t)+(\beta(t)\sin\theta(t)+\dot\lambda(t))\cos\theta(t)\right)J(\nu(t))\\
%&+\left((\beta(t)\cos\theta(t)+\lambda(t)(\dot\theta(t)+\ell(t)))\cos\theta(t)+(\beta(t)\sin\theta(t)+\dot\lambda(t))\sin\theta(t)\right)J(\mu(t))\\
%-\overline{\beta}(t)\overline{\nu}(t)=& \ -\overline{\beta}(t)\left(-\sin\tau(t) J(\overline{\bw}(t))+\cos\tau(t)\overline{\bw}(t)\right)\\
%=& \ \left(-(\beta(t)\cos\theta(t)+\lambda(t)(\dot\theta(t)+\ell(t)))\sin\theta(t)+(\beta(t)\sin\theta(t)+\dot\lambda(t))\cos\theta(t)\right)\mu(t)\\
%&-\left((\beta(t)\cos\theta(t)+\lambda(t)(\dot\theta(t)+\ell(t)))\cos\theta(t)+(\beta(t)\sin\theta(t)+\dot\lambda(t))\sin\theta(t)\right)\nu(t).
%\end{align*}
%By a direct calculation, we have $\overline{\beta}(t)\overline{\bw}(t)=\overline{\beta}(t)(\sin\theta(t)\mu(t)+\cos\theta(t)\nu(t))=\overline{\beta}(t)\bv(t)$. 
By a direct calculation, we have 
\begin{align*}
\overline{\bw}(t)&=\cos\tau(t)\overline{\nu}(t)+\sin\tau(t)\overline{\mu}(t)\\
&=\cos\tau(t)\left(\cos(\theta(t)-\tau(t))\nu(t)+\sin(\theta(t)-\tau(t))\mu(t)\right)\\
&\quad+\sin\tau(t)\left(-\sin(\theta(t)-\tau(t))\nu(t)+\cos(\theta(t)-\tau(t))\mu(t)\right)\\
&=\cos(\tau(t)+\theta(t)-\tau(t))\nu(t)+\sin(\tau(t)+\theta(t)-\tau(t))\mu(t)\\
&=\cos\theta(t)\nu(t)+\sin\theta(t)\mu(t)\\
&=\bv(t).
\end{align*}
It follows that $(\gamma,\nu)$ is a $(\bv,\overline{\bw})$-Bertrand Legendre curve.
\enD
%\begin{align*}
%&\overline{\beta}(t)\overline{\bw}(t)\\
%=& \ \biggl(-\left(-(\beta(t)\cos\theta(t)+\lambda(t)(\dot\theta(t)+\ell(t)))\sin\theta(t)+(\beta(t)\sin\theta(t)+\dot\lambda(t))\cos\theta(t)\right)\cos\tau(t)\\
%& +\left((\beta(t)\cos\theta(t)+\lambda(t)(\dot\theta(t)+\ell(t)))\cos\theta(t)+(\beta(t)\sin\theta(t)+\dot\lambda(t))\sin\theta(t)\right)\sin\tau(t)\biggr)\mu(t)

%\overline{\nu}(t)=& \ \frac{(\beta(t)\cos\theta(t)+\lambda(t)(\dot\theta(t)+\ell(t)))\cos\theta(t)+(\beta(t)\sin\theta(t)+\dot\lambda(t))\sin\theta(t)}{(\beta(t)\cos\theta(t)+\lambda(t)(\dot\theta(t)+\ell(t)))\cos\tau(t)+(\beta(t)\sin\theta(t)+\dot\lambda(t))\sin \tau(t)}\nu(t)\\
%&+\frac{(\beta(t)\cos\theta(t)+\lambda(t)(\dot\theta(t)+\ell(t)))\sin\theta(t)-(\beta(t)\sin\theta(t)+\dot\lambda(t))\cos\theta(t)}{(\beta(t)\cos\theta(t)+\lambda(t)(\dot\theta(t)+\ell(t)))\cos\tau(t)+(\beta(t)\sin\theta(t)+\dot\lambda(t))\sin \tau(t)}\mu(t),\\
%\overline{\mu}(t)=& \ -\frac{(\beta(t)\cos\theta(t)+\lambda(t)(\dot\theta(t)+\ell(t)))\sin\theta(t)-(\beta(t)\sin\theta(t)+\dot\lambda(t))\cos\theta(t)}{(\beta(t)\cos\theta(t)+\lambda(t)(\dot\theta(t)+\ell(t)))\cos\tau(t)+(\beta(t)\sin\theta(t)+\dot\lambda(t))\sin \tau(t)}\nu(t)\\
%&+\frac{(\beta(t)\cos\theta(t)+\lambda(t)(\dot\theta(t)+\ell(t)))\cos\theta(t)+(\beta(t)\sin\theta(t)+\dot\lambda(t))\sin\theta(t)}{(\beta(t)\cos\theta(t)+\lambda(t)(\dot\theta(t)+\ell(t)))\cos\tau(t)+(\beta(t)\sin\theta(t)+\dot\lambda(t))\sin \tau(t)}\mu(t),
%\end{align*}
%It follows that $(\gamma,\nu)$ is a $(\bv,\overline{\bw})$-Bertrand Legendre curve.
%%%%%
%%%%%
\begin{proposition}\label{curvature-Bertrand-Legendre-curve}
Suppose that $(\gamma,\nu)$ and $(\overline{\gamma},\overline{\nu}):  I \to \R^2 \times S^1$ are {$(\bv,\overline{\bw})$-mates}, where $\overline{\gamma}(t)=\gamma(t)+\lambda(t)\bv(t)$, $\bv(t)=\cos \theta(t)\nu(t)+\sin \theta(t) \mu(t)$, $\overline{\bw}(t)=\cos \tau(t) \overline{\nu}(t)+\sin \tau (t) \overline{\mu}(t)$ and $\overline{\nu}(t)=\cos(\theta(t)-\tau(t))\nu(t)+\sin(\theta(t)-\tau(t))\mu(t)$. 
Then the curvature $(\overline{\ell}, \overline{\beta})$ of $(\overline{\gamma},\overline{\nu})$ is given by 
\begin{align*}
\overline{\ell}(t)&=\dot\theta(t)-\dot\tau(t)+\ell(t), \\
\overline{\beta}(t)&=(\beta(t)\cos\theta(t)+\lambda(t)(\dot\theta(t)+\ell(t)))\cos\tau(t)+(\beta(t)\sin\theta(t)+\dot\lambda(t))\sin \tau(t). 
\end{align*}
\end{proposition}
%%%%%
\demo
By differentiating $\overline{\nu}(t)=\cos(\theta(t)-\tau(t))\nu(t)+\sin(\theta(t)-\tau(t))\mu(t)$, we have
$$
\dot{\overline{\nu}}(t)=(\dot\theta(t)-\dot\tau(t)+\ell(t))\left(-\sin(\theta(t)-\tau(t))\nu(t)+\cos(\theta(t)-\tau(t))\mu(t)\right).
$$
It follows that we have 
$\overline{\ell}(t)=\dot{\overline{\nu}}(t)\cdot\overline{\mu}(t)=\dot\theta(t)-\dot\tau(t)+\ell(t)$. 
By Theorem \ref{vw-Legendre}, we have 
$$
\overline{\beta}(t)=\dot{\overline{\gamma}}(t)\cdot\overline{\mu}(t)=(\beta(t)\cos\theta(t)+\lambda(t)(\dot\theta(t)+\ell(t)))\cos\tau(t)+(\beta(t)\sin\theta(t)+\dot\lambda(t))\sin \tau(t).
$$
\enD
%%%%%

%%%%%
\begin{remark}{\rm
Using Propositions \ref{singular-point:criteria} and \ref{curvature-Bertrand-Legendre-curve}, we have criteria for the singular points of $(\overline{\gamma},\overline{\nu})$ with respect to the curvature  $(\overline{\ell}, \overline{\beta})$. 
}
\end{remark}
%%%%%
\par
We give relations between Bertrand Legendre curves and Bertrand regular curves.

%%%%%
\begin{proposition}\label{relations-BL-BR}
$(1)$ Let $(\gamma,\nu)$ and $(\overline{\gamma},\overline{\nu}):I \to \R^2 \times S^1$ be Legendre curves with curvatures $(\ell,\beta)$ and $(\overline{\ell},\overline{\beta})$, respectively. 
Suppose that $(\gamma,\nu)$ and $(\overline{\gamma},\overline{\nu})$ are $(\bv,\overline{\bw})$-mates, where $\bv=\cos \theta \nu+\sin \theta \mu, \overline{\bw}=\cos \tau \overline{\nu}+\sin \tau \overline{\mu}$ and smooth functions $\theta,\tau:I \to \R$. 
Moreover, suppose that $\gamma$ and $\overline{\gamma}$ are regular plane curves.  
Then $\gamma$ and $\overline{\gamma}$ are $(\bv,\overline{\bw})$-regular mates, where 
 $\bv=\sign (\beta) (\cos (\theta-\pi/2) \bt+\sin (\theta-\pi/2) \bn), \overline{\bw}=\sign (\overline{\beta}) (\cos (\tau-\pi/2) \overline{\bt}+\sin (\tau-\pi/2) \overline{\bn})$ and 
$$
\sign(\beta)=\begin{cases}
1 & {\rm if}\ \beta>0 \\
-1 & {\rm if}\ \beta<0 
\end{cases},\quad 
\sign(\overline{\beta})=\begin{cases}
1 & {\rm if}\ \overline{\beta}>0 \\
-1 & {\rm if}\ \overline{\beta}<0 
\end{cases}.
$$
\par
$(2)$ Let $\gamma$ and $\overline{\gamma}:I \to \R^2$ be regular plane curves with curvatures $\kappa$ and $\overline{\kappa}, respectively$. 
Suppose that $\gamma$ and $\overline{\gamma}$ are $(\bv,\overline{\bw})$-regular mates, where $\bv=\cos \theta \bt+\sin \theta \bn, \overline{\bw}=\cos \tau \overline{\bt}+\sin \tau \overline{\bn}$ and smooth functions $\theta,\tau:I \to \R$. 
Then $(\gamma,\bn)$ and $(\overline{\gamma},\overline{\bn})$ are $(\bv,\overline{\bw})$-mates, where  $\bv=\cos (\theta-\pi/2) \bn-\sin (\theta-\pi/2) J(\bn), \overline{\bw}=\cos (\tau-\pi/2) \overline{\bn}-\sin (\tau-\pi/2) J(\overline{\bn})$.
\end{proposition}
%%%%%
\demo
$(1)$ By assumption, $\overline{\gamma}=\gamma+\lambda \bv$ and $\bv=\overline{\bw}$. 
By $\dot{\gamma}=\beta \mu$, we have $|\dot{\gamma}| \bt=\beta \mu$ and $|\dot{\gamma}| \bn=-\beta \nu$. 
Therefore, $|\dot{\gamma}|=|\beta|$, $\mu=\sign (\beta) \bt, \nu=\sign (\beta)\bn$. 
It follows that 
\begin{align*}
\bv &=\cos \theta \nu +\sin \theta \mu=\sign (\beta)(-\cos \theta \bn+\sin \theta \bt)\\
&= \sign (\beta) \left(\cos (\theta-\pi/2)\bt+\sin (\theta-\pi/2)\bn \right).
\end{align*}
We can also calculate 
$$
\overline{\bw} =\cos \tau \overline{\nu} +\sin \tau \overline{\mu}= \sign (\overline{\beta}) \left(\cos (\tau-\pi/2) \overline{\bt}+\sin (\tau-\pi/2)\overline{\bn} \right).
$$
Then $\gamma$ and $\overline{\gamma}$ are $(\bv,\overline{\bw})$-regular curves. 
\par
$(2)$ By assumption, $\overline{\gamma}=\gamma+\lambda \bv$ and $\bv=\overline{\bw}$. 
Since 
\begin{align*}
\bv &= \cos \theta \bt + \sin \theta \bn=\cos (\theta-\pi/2)\bn -\sin (\theta-\pi/2) J(\bn),\\
\overline{\bw} &= \cos \tau \overline{\bt} + \sin \tau \overline{\bn}=\cos (\tau-\pi/2) \overline{\bn}-\sin (\tau-\pi/2) J(\overline{\bn}), 
\end{align*}
$(\gamma,\bn)$ and $(\overline{\gamma},\overline{\bn})$ are $(\bv,\overline{\bw})$-mates.
\enD
%%%%%
By Theorem \ref{vw-Legendre}, we have the special cases of Bertrand Legendre curves.
%%%%%
\begin{corollary}\label{vw-Legendre-special}
Under the same notations as in Theorem \ref{vw-Legendre}, we have the following.
\par
$(1)$ Suppose that $\theta(t)=0$ and $\tau(t)=0$ for all $t \in I$. 
Then $(\gamma,\nu)$ is a $(\nu,\overline{\nu})$-Bertrand Legendre curve if and only if $\lambda(t)$ is a constant. 
It follows that $\overline{\gamma}$ is a parallel curve of $(\gamma,\nu)$.
\par
$(2)$ Suppose that $\theta(t)=0$ and $\tau(t)=\pi/2$ for all $t \in I$. 
Then $(\gamma,\nu)$ is a $(\nu,\overline{\mu})$-Bertrand Legendre curve if and only if $\beta(t)+\lambda(t) \ell(t)=0$ for all $t \in I$. 
It follows that $\overline{\gamma}$ is an evolute of $(\gamma,\nu)$.
\par
$(3)$ Suppose that $\theta(t)=\pi/2$ and $\tau(t)=0$ for all $t \in I$. 
Then $(\gamma,\nu)$ is a $(\mu,\overline{\nu})$-Bertrand Legendre curve if and only if $\beta(t)+\dot{\lambda}(t)=0$ for all $t \in I$. 
It follows that $\overline{\gamma}$ is an involute of $(\gamma,\nu)$.
\par
$(4)$ Suppose that $\theta(t)=\pi/2$ and $\tau(t)=\pi/2$ for all $t \in I$. 
Then $(\gamma,\nu)$ is a $(\mu,\overline{\mu})$-Bertrand Legendre curve if and only if $\ell(t)=0$ for all $t \in I$. 
It follows that $\overline{\gamma}$ and $\gamma$ are a part of line.
\par
$(5)$ Suppose that $\theta(t)$ is a constant $\theta$ and $\tau(t)=\pi/2$ for all $t \in I$. 
Then $(\gamma,\nu)$ is a $(\bv,\overline{\mu})$-Bertrand Legendre curve if and only if 
$\beta(t) \cos \theta+\lambda(t) \ell(t)=0$ for all $t \in I$. 
\par
$(6)$ Suppose that $\theta(t)=\pi/2$ and $\tau(t)$ is a constant $\tau$ for all $t \in I$. 
Then $(\gamma,\nu)$ is a $(\mu,\overline{\bw})$-Bertrand Legendre curve if and only if $(\beta(t)+\dot{\lambda}(t)) \cos \tau -\lambda(t)\ell(t)\sin \tau=0$ for all $t \in I$. 
\end{corollary}
%%%%%

Using Corollary \ref{vw-Legendre-special}, we may directly define the evolutoid and involutoid of Legendre curves. 
%%%%%
\begin{definition}\label{evolutoid-involutoid}{\rm 
Let $(\gamma,\nu):I \to \R^2 \times S^1$ be a Legendre curve and $\theta, \tau$ be constants. 
We say that 
$$
\mathcal{E}v[\theta](\gamma,\nu):I \to \R^2, \ \mathcal{E}v[\theta](\gamma,\nu)(t)=\gamma(t)+\lambda(t)(\cos \theta \nu(t)+\sin \theta \mu(t)),
$$ 
where $\beta(t) \cos \theta+\lambda(t) \ell(t)=0$ for all $t \in I$ is an {\it evolutoid} ({\it $\theta$-evolutoid}) of the Legendre curve $(\gamma,\nu)$ and 
$$
\mathcal{I}nv[\tau](\gamma,\nu):I \to \R^2, \ \mathcal{I}nv[\tau](\gamma,\nu)(t)=\gamma(t)+\lambda(t)\mu(t),
$$ where 
$(\beta(t)+\dot{\lambda}(t)) \cos \tau -\lambda(t)\ell(t)\sin \tau=0$ for all $t \in I$ is an {\it involutoid} ({\it $\tau$-involutoid}) of the Legendre curve $(\gamma,\nu)$.
}
\end{definition}
%%%%%

Note that the definitions of the evolutoid of regular plane curves and frontals have already investigated in \cite{Giblin-Warder, Izumiya-Takeuchi-2019}, and of the involutoid (tanvolute) of regular plane curves have already investigated in \cite{AM}. 
Moreover, the evolutoid and involutoid of spherical Legendre curves investigated in \cite{Li-Pei}.
However, the explicit form of the definition of the involutoid of regular plane curves and of Legendre curves (frontals)  in the unit tangent bundle over Euclidean plane can not find  as far as we know. 

%%%%%
\begin{corollary}
Let $(\gamma,\nu):I \to \R^2 \times S^1$ be a Legendre curve with curvature $(\ell,\beta)$ and $\theta, \tau$ be constants.
\par
$(1)$ $(\mathcal{E}v[\theta], \nu_{\mathcal{E}v}[\theta])(\gamma,\nu):I \to \R^2 \times S^1$ is a Legendre curve with curvature $(\ell_{\mathcal{E}v}[\theta],\beta_{\mathcal{E}v}[\theta])=(\ell,\beta \sin \theta+\dot{\lambda})$, where 
$\nu_{\mathcal{E}v}[\theta](\gamma,\nu)=\sin \theta \nu-\cos \theta \mu$ and $\beta \cos \theta+\lambda \ell=0$.
\par
$(2)$ $(\mathcal{I}nv[\tau], \nu_{\mathcal{I}nv}[\tau])(\gamma,\nu):I \to \R^2 \times S^1$ is a Legendre curve with curvature $(\ell_{\mathcal{I}nv}[\tau],\beta_{\mathcal{I}nv}[\tau])=(\ell,\lambda\ell\cos \tau+(\beta +\dot{\lambda})\sin \tau)$, where 
$\nu_{\mathcal{I}nv}[\tau](\gamma,\nu)=\sin \tau \nu+\cos \tau \mu$ and $(\beta+\dot{\lambda}) \cos \tau -\lambda \ell \sin \tau=0$.
\end{corollary}
%%%%%

%%%%%
\begin{remark}{\rm
We can observe that the inflection points of $(\gamma,\nu)$ are invariants for evolutoids and involutoids. 
That is, if $t_0$ is an inflection point of $(\gamma,\nu)$, then 
$t_0$ is also an inflection point of $(\mathcal{E}v[\theta], \nu_{\mathcal{E}v}[\theta])(\gamma,\nu)$ and $(\mathcal{I}nv[\tau], \nu_{\mathcal{I}nv}[\tau])(\gamma,\nu)$.
On the other hand, the singular points of $\mathcal{E}v[\theta](\gamma,\nu)$ and $\mathcal{I}nv[\tau](\gamma,\nu)$ may be changed for singular points of $\gamma$. 
Moreover,  $\mathcal{E}v[\theta](\gamma,\nu)$ connects between evolutes $(\theta=0)$ and the original curve $\gamma$ $(\theta=\pi/2)$ if $\ell \not\equiv 0$, and $\mathcal{I}nv[\tau](\gamma,\nu)$ connects between involutes $(\tau=0)$ and the original curve $(\tau=\pi/2)$ if $\ell \not\equiv 0$.
}
\end{remark}
%%%%%
\par
We give new correspondences which connects between evolutes and involutes of Legendre curves. 
That is, we consider $(\bv,J(\overline{\bv}))$ and $(J(\bw),\overline{\bw})$-Bertrand Legendre curves.

%%%%%
\begin{definition}\label{nvolute}{\rm 
Let $(\gamma,\nu):I \to \R^2 \times S^1$ be a Legendre curve with curvature $(\ell,\beta)$, and $\theta, \tau$ be constants. 
\par
$(1)$ We define
$$
N[\theta](\gamma,\nu):I \to \R^2, \ N[\theta](\gamma,\nu)(t)=\gamma(t)+\lambda(t)(\cos \theta \nu(t)+\sin \theta \mu(t)),
$$
where $\beta(t)+\dot{\lambda}(t)\sin \theta+\lambda(t)\ell(t)\cos \theta=0$ for all $t \in I$. 
\par
$(2)$ We define
$$
T[\tau](\gamma,\nu):I \to \R^2, \ T[\tau](\gamma,\nu)(t)=\gamma(t)+\lambda(t)(\cos \tau \mu(t)-\sin \tau \nu(t)),
$$
where $\beta(t)+\dot{\lambda}(t) \cos \tau-\lambda(t) \ell(t) \sin \tau=0$ for all $t \in I$.
}
\end{definition}
%%%%%
\par
If $\theta=0$, then $N[0](\gamma,\nu)$ is an evolute of $(\gamma,\nu)$ and if $\theta=\pi/2$, then $N[\pi/2](\gamma,\nu)$ is an involute of $(\gamma,\nu)$.
Moreover, if $\tau=0$, then $T[0](\gamma,\nu)$ is an involute of $(\gamma,\nu)$ and if $\tau=-\pi/2$, then $T[-\pi/2](\gamma,\nu)$ is an evolute of $(\gamma,\nu)$.

%%%%%
\begin{corollary}\label{nvolute-curvature}
Let $(\gamma,\nu):I \to \R^2 \times S^1$ be a Legendre curve with curvature $(\ell,\beta)$ and $\theta, \tau$ be constants.
\par
$(1)$ $(N[\theta], \nu_{N}[\theta])(\gamma,\nu):I \to \R^2 \times S^1$ is a Legendre curve with curvature $(\ell_{N}[\theta],\beta_{N}[\theta])=(\ell, -\lambda \ell \sin \theta+\dot{\lambda} \cos \theta)$, where 
$\nu_{N}[\theta](\gamma,\nu)=-\mu$ and $\beta +\dot{\lambda} \sin \theta+\lambda \ell \cos \theta=0$.
\par
$(2)$ $(T[\tau]),\nu_{T}[\tau](\gamma,\nu):I \to \R^2 \times S^1$ is a 
Legendre curve with curvature $(\ell_{T}[\tau],\beta_{T}[\tau])=(\ell, \lambda \ell \cos \tau+\dot{\lambda} \sin \tau)$, where 
$\nu_{T}[\tau](\gamma,\nu)=\mu$ and $\beta +\dot{\lambda} \cos \tau-\lambda \ell \sin \tau=0$. 
\end{corollary}
%%%%%

%%%%%
We give an inverse operation of Bertrand Legendre curves. 
The set of Legendre curves is denoted by 
$$
L(I,\R^2 \times S^1) :=\{ (\gamma,\nu) \in C^{\infty}(I, \R^2 \times S^1)\ |\ \dot{\gamma}(t) \cdot \nu(t)=0 \ {\rm for \ all} \ t \in I\}.
$$
\par
We consider an operator between Legendre curves, 
$$
L[\bv,\overline{\bw}]=L[\theta,\tau]: L(I,\R^2 \times S^1) \to L(I,\R^2 \times S^1)
$$ by 
$$
L[\bv,\overline{\bw}](\gamma,\nu)=L[\theta,\tau](\gamma,\nu)=(\gamma+\lambda(\cos \theta \nu +\sin \theta \mu), \cos (\theta-\tau) \nu+\sin (\theta-\tau)\mu),
$$
where $\lambda:I \to \R$ is satisfied 
$$
(\beta(t) \sin \theta(t)+\dot{\lambda}(t)) \cos \tau(t)-(\beta(t)\cos \theta(t)+\lambda(t)(\dot{\theta}(t)+\ell(t)))\sin \tau(t)=0
$$
for all $t \in I$. 
By Proposition \ref{curvature-Bertrand-Legendre-curve}, the curvature $(\ell[\bv,\overline{\bw}], \beta[\bv,\overline{\bw}])=(\ell[\theta,\tau],\beta[\theta,\tau])$ of the Legendre curve $L[\bv,\overline{\bw}]=L[\theta,\tau]$ is given by 
\begin{align*}
\ell[\bv,\overline{\bw}](t)&=\dot\theta(t)-\dot\tau(t)+\ell(t), \\
\beta[\bv,\overline{\bw}](t)&=(\beta(t)\cos\theta(t)+\lambda(t)(\dot\theta(t)+\ell(t)))\cos\tau(t)+(\beta(t)\sin\theta(t)+\dot\lambda(t))\sin \tau(t). 
\end{align*}

%%%%%
\begin{lemma}\label{relation}
Let $(\gamma_i,\nu_i):I \to \R^2 \times S^1$ be Legendre curves for $i=1,2,3$.
Suppose that $(\gamma_1,\nu_1)$ and $(\gamma_2,\nu_2)$ are $(\bv_1,\bv_2)$-mates with $\lambda_1$, and $(\gamma_2,\nu_2)$ and $(\gamma_3,\nu_3)$ are $(\bv_2,\bv_3)$-mates with $\lambda_2$.
\par
$(1)$ If $\lambda_1(t)+\lambda_2(t)=0$ for all $t \in I$, then $\gamma_1=\gamma_3$.
\par
$(2)$ If $\lambda_1+\lambda_2 \not\equiv 0$, then $(\gamma_1,\nu_1)$ and $(\gamma_3,\nu_3)$ are $(\bv_1,\bv_3)$-mates.
\end{lemma}
%%%%%
\demo
By the assumptions, we have $\gamma_2(t)=\gamma_1(t)+\lambda_1(t)\bv_1(t), \bv_1(t)=\bv_2(t)$ and $\gamma_3(t)=\gamma_2(t)+\lambda_2(t)\bv_2(t), \bv_2(t)=\bv_3(t)$ for all $t \in I$. 
Then 
\begin{align*}
\gamma_3(t) &=\gamma_1(t) +\lambda_1(t) \bv_1(t)+\lambda_2(t)\bv_2(t)=\gamma_1(t)+(\lambda_1(t)+\lambda_2(t))\bv_1(t), \\
\bv_1(t) &=\bv_3(t)
\end{align*}
for all $t \in I$. 
It follows that if $\lambda_1(t)+\lambda_2(t)=0$ for all $t \in I$, then $\gamma_1(t)=\gamma_3(t)$ for all $t \in I$. 
\par
Moreover, by Definition \ref{Bertrand-type-Legendre}, if $\lambda_1+\lambda_2 \not\equiv 0$, then $(\gamma_1,\nu_1)$ and $(\gamma_3,\nu_3)$ are $(\bv_1,\bv_3)$-mates.
\enD
%%%%%

%%%%%
\begin{theorem}\label{inverse-relation} 
Let $(\gamma,\nu)$ and $(\overline{\gamma},\overline{\nu}):I \to \R^2 \times S^1$ be Legendre curves. 
Suppose that $(\gamma,\nu)$ and $(\overline{\gamma},\overline{\nu})$ are $(\bv,\overline{\bw})$-mates with $\lambda$. 
We denote $(\overline{\gamma},\overline{\nu})=L[\bv,\overline{\bw}](\gamma,\nu)=L[\theta,\tau](\gamma,\nu)$. 
\par
$(1)$ $(\overline{\gamma},\overline{\nu})$ and $(\gamma,\nu)$ are $(\overline{\bw},\bv)$-mates with $-\lambda$. 
Moreover, $L[\bw,\overline{\bv}] \circ L[\bv,\overline{\bw}](\gamma,\nu)=L[\tau,\theta] \circ L[\theta,\tau](\gamma,\nu)=(\gamma,\nu)$.
\par
$(2)$ Suppose that $(\overline{\gamma},\overline{\nu})$ and $(\gamma,\nu)$ are $(\overline{\bw},\bv)$-mates with $\overline{\lambda}$. 
If $\lambda +\overline{\lambda} \not\equiv 0$, then $(\gamma,\nu)$ and $L[\bw,\overline{\bv}] \circ L[\bv,\overline{\bw}](\gamma,\nu)=L[\tau,\theta] \circ L[\theta,\tau](\gamma,\nu)$ are $(\bv,\overline{\bv})$-mates.  
\end{theorem}
%%%%%
\demo 
$(1)$ By assumption, we have $\overline{\gamma}=\gamma+\lambda \bv$ and $\bv=\overline{\bw}$. 
It follows that $\gamma=\overline{\gamma}-\lambda \overline{\bw}$ and $\overline{\bw}=\bv$.
Hence, $(\overline{\gamma},\overline{\nu})$ and $(\gamma,\nu)$ are $(\overline{\bw},\bv)$-mates with $-\lambda$. 
Moreover, 
\begin{align*}
L[\bw,\overline{\bv}] \circ L[\bv,\overline{\bw}](\gamma,\nu) &=L[\bw,\overline{\bv}](\overline{\gamma},\overline{\nu})\\
&=L[\bw,\overline{\bv}](\gamma+\lambda \bv,\cos (\theta-\tau)\nu+\sin (\theta-\tau)\mu)\\
&=(\overline{\gamma}-\lambda \overline{\bw}, \cos (\tau-\theta)\overline{\nu}+\sin (\tau-\theta)\overline{\mu}).
\end{align*} 
By Lemma \ref{relation} $(1)$, we have $\overline{\gamma}-\lambda \overline{\bw}=\gamma$. 
By a direct calculation, 
\begin{align*}
\cos (\tau-\theta)\overline{\nu}+\sin (\tau-\theta)\overline{\mu}
&= \cos (\tau-\theta)(\cos (\theta-\tau)\nu+\sin (\theta-\tau)\mu)\\
&\quad +\sin (\tau-\theta)(\cos (\theta-\tau)\mu-\sin (\theta-\tau)\nu)\\
&=\nu.
\end{align*}
\par
$(2)$ By Lemma \ref{relation} $(2)$, we have the result. 
\enD
%%%%%
By Theorem \ref{inverse-relation}, we have the following Corollary.
%%%%%
\begin{corollary}\label{evolutoid-involutoid-inverse} 
Let $(\gamma,\nu)$ and $(\overline{\gamma},\overline{\nu}):I \to \R^2 \times S^1$ be Legendre curves and $\theta$ be a constant.  
\par
$(1)$ Suppose that $(\gamma,\nu)$ and $(\overline{\gamma},\overline{\nu})$ are $(\bv,\overline{\mu})$-mates with $\lambda$. 
We denote $(\overline{\gamma},\overline{\nu})=L[\bv,\overline{\mu}](\gamma,\nu)=L[\theta,\pi/2](\gamma,\nu)=(\mathcal{E}v[\theta], \nu_{\mathcal{E}v}[\theta])(\gamma,\nu)$. 
\par
$(i)$ $(\overline{\gamma},\overline{\nu})$ and $(\gamma,\nu)$ are $(\overline{\mu},\bv)$-mates with $-\lambda$. 
Moreover, 
\begin{align*}
L[\mu,\overline{\bv}] \circ L[\bv,\overline{\mu}](\gamma,\nu) &=L\left[\frac{\pi}{2},\theta\right] \circ L\left[\theta,\frac{\pi}{2}\right](\gamma,\nu) \\
&=(\mathcal{I}nv[\theta],\nu_{\mathcal{I}nv}[\theta]) \circ (\mathcal{E}v[\theta],\nu_{\mathcal{E}v}[\theta])(\gamma,\nu)=(\gamma,\nu).
\end{align*}
\par
$(ii)$ Suppose that $(\overline{\gamma},\overline{\nu})$ and $(\gamma,\nu)$ are $(\overline{\mu},\bv)$-mates with $\overline{\lambda}$. 
If $\lambda +\overline{\lambda} \not\equiv 0$, then $(\gamma,\nu)$ and 
$$
L[\mu,\overline{\bv}] \circ L[\bv,\overline{\mu}](\gamma,\nu)=L\left[\frac{\pi}{2},\theta\right] \circ L\left[\theta,\frac{\pi}{2}\right](\gamma,\nu)=(\mathcal{I}nv[\theta],\nu_{\mathcal{I}nv}[\theta]) \circ (\mathcal{E}v[\theta],\nu_{\mathcal{E}v}[\theta])(\gamma,\nu)
$$ 
are $(\bv,\overline{\bv})$-mates.  
\par
$(2)$ Suppose that $(\gamma,\nu)$ and $(\overline{\gamma},\overline{\nu})$ are $(\mu,\overline{\bv})$-mates with $\lambda$. 
We denote $(\overline{\gamma},\overline{\nu})=L[\mu,\overline{\bv}](\gamma,\nu)=L[\pi/2,\theta](\gamma,\nu)=(\mathcal{I}nv[\theta], \nu_{\mathcal{I}nv}[\theta])(\gamma,\nu)$. 
\par
$(i)$ $(\overline{\gamma},\overline{\nu})$ and $(\gamma,\nu)$ are $(\overline{\bv},\mu)$-mates with $-\lambda$. 
Moreover, 
\begin{align*}
L[\bv,\overline{\mu}] \circ L[\mu,\overline{\bv}](\gamma,\nu) &=L\left[\theta,\frac{\pi}{2}\right] \circ L\left[\frac{\pi}{2},\theta\right](\gamma,\nu) \\
&=(\mathcal{E}v[\theta],\nu_{\mathcal{E}v}[\theta]) \circ (\mathcal{I}nv[\theta],\nu_{\mathcal{I}nv}[\theta]) (\gamma,\nu)=(\gamma,\nu).
\end{align*}
\par
$(ii)$ Suppose that $(\overline{\gamma},\overline{\nu})$ and $(\gamma,\nu)$ are $(\overline{\bv},\mu)$-mates with $\overline{\lambda}$. 
If $\lambda +\overline{\lambda} \not\equiv 0$, then $(\gamma,\nu)$ and 
\begin{align*}
L[\bv,\overline{\mu}] \circ L[\mu,\overline{\bv}](\gamma,\nu) =L\left[\theta,\frac{\pi}{2}\right] \circ L\left[\frac{\pi}{2},\theta\right](\gamma,\nu) =(\mathcal{E}v[\theta],\nu_{\mathcal{E}v}[\theta]) \circ (\mathcal{I}nv[\theta],\nu_{\mathcal{I}nv}[\theta]) (\gamma,\nu)
\end{align*} are $(\mu,\overline{\mu})$-mates. 
\par
$(3)$ Suppose that $(\gamma,\nu)$ and $(\overline{\gamma},\overline{\nu})$ are $(\bv,J(\overline{\bv}))$-mates with $\lambda$. 
We denote $(\overline{\gamma},\overline{\nu})=L[\bv,J(\overline{\bv})](\gamma,\nu)=L[\theta,\theta+\pi/2](\gamma,\nu)=(N[\theta], \nu_{N}[\theta])(\gamma,\nu)$. 
\par
$(i)$ $(\overline{\gamma},\overline{\nu})$ and $(\gamma,\nu)$ are $(J(\overline{\bv}),\bv)$-mates with $-\lambda$. 
Moreover, 
\begin{align*}
L[J(\bv),\overline{\bv}] \circ L[\bv,J(\overline{\bv})](\gamma,\nu) &=L\left[\theta+\frac{\pi}{2},\theta\right] \circ L\left[\theta,\theta+\frac{\pi}{2}\right](\gamma,\nu) \\
&=(T[\theta],\nu_{T}[\theta]) \circ (N[\theta],\nu_{N}[\theta])(\gamma,\nu)=(\gamma,\nu).
\end{align*}
\par
$(ii)$ Suppose that $(\overline{\gamma},\overline{\nu})$ and $(\gamma,\nu)$ are $(J(\overline{\bv}),\bv)$-mates with $\overline{\lambda}$. 
If $\lambda +\overline{\lambda} \not\equiv 0$, then $(\gamma,\nu)$ and 
$$
L[J(\bv),\overline{\bv}] \circ L[\bv,J(\overline{\bv})](\gamma,\nu)=L\left[\theta+\frac{\pi}{2},\theta\right] \circ L\left[\theta,\theta+\frac{\pi}{2}\right](\gamma,\nu)=(T[\theta],\nu_{T}[\theta]) \circ (N[\theta],\nu_{N}[\theta])(\gamma,\nu)
$$ are $(\bv,\overline{\bv})$-mates.  
\par
$(4)$ Suppose that $(\gamma,\nu)$ and $(\overline{\gamma},\overline{\nu})$ are $(J(\bv),\overline{\bv})$-mates with $\lambda$. 
We denote $(\overline{\gamma},\overline{\nu})=L[J(\bv),\overline{\bv}](\gamma,\nu)=L[\theta+\pi/2,\theta](\gamma,\nu)=(T[\theta], \nu_{T}[\theta])(\gamma,\nu)$. 
\par
$(i)$ $(\overline{\gamma},\overline{\nu})$ and $(\gamma,\nu)$ are $(\overline{\bv}),J(\bv))$-mates with $-\lambda$. 
Moreover, 
\begin{align*}
L[\bv,J(\overline{\bv})] \circ L[J(\bv),\overline{\bv}](\gamma,\nu) &=L\left[\theta,\theta+\frac{\pi}{2}\right] \circ L\left[\theta+\frac{\pi}{2},\theta\right](\gamma,\nu) \\
&=(N[\theta],\nu_{N}[\theta]) \circ (T[\theta],\nu_{T}[\theta])(\gamma,\nu)=(\gamma,\nu).
\end{align*}
\par
$(ii)$ Suppose that $(\overline{\gamma},\overline{\nu})$ and $(\gamma,\nu)$ are $(\overline{\bv},J(\bv))$-mates with $\overline{\lambda}$. 
If $\lambda +\overline{\lambda} \not\equiv 0$, then $(\gamma,\nu)$ and 
$$
L[\bv,J(\overline{\bv})] \circ L[J(\bv),\overline{\bv}](\gamma,\nu)=L\left[\theta,\theta+\frac{\pi}{2}\right] \circ L\left[\theta+\frac{\pi}{2},\theta\right](\gamma,\nu)=(N[\theta],\nu_{N}[\theta]) \circ (T[\theta],\nu_{T}[\theta])(\gamma,\nu)
$$ are $(J(\bv),J(\overline{\bv}))$-mates.  
\end{corollary}
%%%%%

We consider a subset of the set of Legendre curves.
Let $\bv,\overline{\bw}:I \to S^1$ be mappings with $\theta,\tau:I \to \R$. 
We denote 
$$
L_S[\bv,\overline{\bw}]=L_S[\theta,\tau]=\{(\gamma,\nu) \in L(I,\R^2 \times S^1) \ | \ {\rm there \ exists} \ \lambda:I \to \R {\rm \ such \ that}\ \eqref{vw-Legendre-condition} \ {\rm holds}\}.
$$

For a mapping $\bu:I \to S^1$, we consider an equivalence relation of $L[\bv,\overline{\bw}]$ (respectively, $L_S[\bv,\overline{\bw}]$). 
Let $(\gamma,\nu)$ and $(\widetilde{\gamma},\widetilde{\nu}) \in L[\bv,\overline{\bw}]$ (respectively, $L_S[\bv,\overline{\bw}]$). 
We define a relation $(\gamma,\nu) \sim (\widetilde{\gamma},\widetilde{\nu})$ if $(\gamma,\nu)$ and $ (\widetilde{\gamma},\widetilde{\nu})$ are $(\bu,\overline{\bu})$-mates with $\lambda$. 
Here we drop the condition $\lambda \not\equiv 0$, that is, we admit the case of $\lambda=0$. 
By Lemma \ref{relation} and Theorem \ref{inverse-relation}, the relation is an equivalence relation.
We consider the quotient space of the set of Legendre curves $L[\bv,\overline{\bw}]$ (respectively, $L_S[\bv,\overline{\bw}]$) by the equivalence relation and denote it by $L[\bv,\overline{\bw}]/\sim [\bu,\overline{\bu}]$ (respectively, $L_S[\bv,\overline{\bw}]/\sim [\bu,\overline{\bu}]$). 
Then we define a mapping between $L_S[\bv,\overline{\bw}]$ and prove that it is bijective up to equivalence relations. 

%%%%%
\begin{theorem}\label{bijective}
Let $(\gamma,\nu):I \to \R^2 \times S^1$ be a Legendre curve with curvature $(\ell,\beta)$ and $\bv,\overline{\bw}:I \to S^1$ be mappings with $\theta,\tau$.
\par
$(1)$ $[L[\bv,\overline{\bw}]]:L_S[\bv,\overline{\bw}] \to \left(L_S[\bw,\overline{\bv}]/\sim [\bw,\overline{\bw}]\right)$, $(\gamma,\nu) \mapsto [L[\bv,\overline{\bw}]](\gamma,\nu):=[L[\bv,\overline{\bw}](\gamma,\nu)]$ is a mapping. 
\par
$(2)$ $[L[\bv,\overline{\bw}]]: \left(L_S[\bv,\overline{\bw}]/\sim [\bv,\overline{\bv}]\right) \to \left(L_S[\bw,\overline{\bv}]/\sim [\bw,\overline{\bw}]\right)$, $[(\gamma,\nu)] \mapsto [L[\bv,\overline{\bw}]]([(\gamma,\nu)]):=[L[\bv,\overline{\bw}](\gamma,\nu)]$ is bijective. 
\end{theorem}
%%%%%
\demo
$(1)$ First, we show that $[L[\bv,\overline{\bw}]](\gamma,\nu) \in L_S[\bw,\overline{\bv}]$ for all $(\gamma,\nu) \in L_S[\bv,\overline{\bw}]$. 
Since $L[\bv,\overline{\bw}](\gamma,\nu)=(\gamma+\lambda (\cos \theta \nu+\sin \theta \mu),\cos (\theta-\tau)\nu+\sin (\theta-\tau)\mu)$ and the curvature is given by $(\overline{\ell},\overline{\beta})=(\ell[\bv,\overline{\bw}],\beta[\bv,\overline{\bw}])$ by Proposition \ref{curvature-Bertrand-Legendre-curve}, 
$\overline{\lambda}=-\lambda$ satisfies the condition 
$$
(\overline{\beta}(t) \sin \tau (t)+\dot{\overline{\lambda}}(t))\cos \theta(t) - 
(\overline{\beta}(t) \cos \tau (t)+\overline{\lambda}(t)(\dot{\tau}(t)+\overline{\ell}(t)))\sin \theta(t)=0
$$
for all $t \in I$.
It follows that $L[\bv,\overline{\bw}](\gamma,\nu) \in L_S[\bw,\overline{\bv}]$.
\par
Next, we show that for all $(\gamma,\nu) \in L_S[\bv,\overline{\bw}]$, there exists unique 
$$
[L[\bv,\overline{\bw}](\gamma,\nu)] \in \left(L_S[\bw,\overline{\bv}]/\sim [\bw,\overline{\bw}]\right).
$$
If there exist $\lambda$ and $\widetilde{\lambda}$ such that condition \eqref{vw-Legendre-condition} satisfy for all $t \in I$. 
Then if we take $\overline{\lambda}=\widetilde{\lambda}-\lambda$, then $L[\bv,\overline{\bw}](\gamma,\nu)$ with $\lambda$ and $L[\bv,\overline{\bw}](\gamma,\nu)$ with $\widetilde{\lambda}$ are $(\bw,\overline{\bw})$-mates with $\overline{\lambda}$. 
It follows that $[L[\bv,\overline{\bw}]]$ is a mapping. 
\par
$(2)$ First, we show that the mapping is well-defined.
Suppose that $(\gamma,\nu)$ and $(\widetilde{\gamma},\widetilde{\nu})$ are $(\bv,\overline{\bv})$-mates with $\widetilde{\lambda}$. 
Then we can show that $L[\bv,\overline{\bw}](\gamma,\nu)$ with $\lambda$ and $L[\bv,\overline{\bw}](\widetilde{\gamma},\widetilde{\nu})$ with $\lambda-\widetilde{\lambda}$ are the same, that is, $L[\bv,\overline{\bw}](\gamma,\nu)=L[\bv,\overline{\bw}](\widetilde{\gamma},\widetilde{\nu})$. 
Therefore, $L[\bv,\overline{\bw}](\gamma,\nu)$ and $L[\bv,\overline{\bw}](\widetilde{\gamma},\widetilde{\nu})$ are $(\bw,\widetilde{\bw})$-mates by $(1)$. 
\par
Next, we show that the mapping is injective. 
Suppose that $L[\bv,\widetilde{\bw}](\gamma,\nu)$ and $L[\bv,\widetilde{\bw}](\widetilde{\gamma},\widetilde{\nu})$ are $(\bw,\overline{\bw})$-mates. 
By Theorem \ref{inverse-relation}, if we consider $L[\bw,\overline{\bv}]$, then we have 
$$
(\gamma,\nu)=L[\bw,\overline{\bv}] \circ L[\bv,\overline{\bw}](\gamma,\nu)
=L[\bw,\overline{\bv}] \circ L[\bv,\overline{\bw}](\widetilde{\gamma},\widetilde{\nu}).
$$
Since $(\widetilde{\gamma},\widetilde{\nu})$ and $L[\bw,\overline{\bv}] \circ L[\bv,\overline{\bw}](\widetilde{\gamma},\widetilde{\nu})$ are $(\bv,\overline{\bv})$-mates, $(\gamma,\nu)$ and $(\widetilde{\gamma},\widetilde{\nu})$ are also $(\bv,\overline{\bv})$-mates. 
\par
Finally, we show that the mapping is surjective. 
For any $[(\overline{\gamma},\overline{\nu})] \in L_S[\bw,\overline{\bv}]$, there exists $\lambda$ with condition \eqref{vw-Legendre-condition} satisfies for all $t \in I$.
We define $(\gamma,\nu) \in L_S[\bv,\overline{\bw}]$ by 
$\gamma=\overline{\gamma}+\overline{\lambda} \overline{\bw}$, 
$\nu=\cos (\tau-\theta) \overline{\nu}+\sin (\tau-\theta) \overline{\mu}$, where $\overline{\lambda}=-\lambda$. 
Then we have $L[\bv,\overline{\bw}](\gamma,\nu)=(\overline{\gamma},\overline{\nu})$. 
It follows that $[L[\bv,\overline{\bw}]]([\gamma,\nu])=[(\overline{\gamma},\overline{\nu})]$.
\enD
%%%%%

Finally, we give concrete examples of Bertrand Legendre curves.
%%%%%
\begin{example}{\rm
Let $(\gamma,\nu):[0,2\pi) \to \R^2 \times S^1$ be 
\begin{align*}
\gamma(t)=\left(r\cos t,r\sin t\right),\ \nu(t)=\left(\cos t,\sin t\right),
\end{align*}
where $r$ is a positive constant.
Then $\gamma$ is a circle. 
By a direct calculation, $\mu(t)=(-\sin t,\cos t)$ and $(\gamma,\nu)$ is a Legendre curve with the curvature $(\ell(t),\beta(t))=(1, r)$. 
\par
If we take  $\theta(t)=0$, $\tau(t)=\pi/2$ and a smooth function $\lambda: [0,2\pi) \to \R$ by $\lambda(t)=-r$, then $\beta(t)+\lambda(t)\ell(t)=0$ for all $t \in [0,2\pi)$. 
Therefore, the Legendre curve $(\overline{\gamma},\overline{\nu})$ such that $\overline{\gamma}$ is an evolute of $(\gamma,\nu)$ is given by
\begin{align*}
\overline{\gamma}(t)=\gamma(t)+\lambda(t)\bv(t)=\left(0,0\right), \ 
\overline{\nu}(t)=-\mu(t)=\left(\sin t,-\cos t\right).
\end{align*}
\par
If we take $\theta(t)=\pi/2$, $\tau(t)=0$ and a smooth function $\lambda: [0,2\pi) \to \R$ by $\lambda(t)=-rt+c$, where $c$ is a constant, then $\beta(t)+\dot{\lambda}(t)=0$ for all $t \in [0,2\pi)$.  
Therefore, the Legendre curve $(\overline{\gamma},\overline{\nu})$ such that $\overline{\gamma}$ is an involute of $(\gamma,\nu)$ is given by
\begin{align*}
\overline{\gamma}(t)&=\gamma(t)+\lambda(t)\bv(t)=\left(r(\cos t+t\sin t)-c\sin t,r(\sin t-t\cos t)+c\cos t\right), \\
\overline{\nu}(t)&=\mu(t)=\left(-\sin t,\cos t\right).
\end{align*}
\par 
If we take $\tau(t)=\pi/2$, a constant $\theta$ and a smooth function $\lambda: [0,2\pi) \to \R$ by $\lambda(t)=-r\cos\theta$,
then $\beta(t) \cos \theta+\lambda(t) \ell(t)=0$ for all $t \in [0,2\pi)$. 
By Definition \ref{evolutoid-involutoid}, the evolutoid $\mathcal{E}v[\theta](\gamma,\nu)$ of $(\gamma,\nu)$ is given by 
\begin{align*}
\mathcal{E}v[\theta](\gamma,\nu)(t)&=\gamma(t)+\lambda(t)(\cos\theta\nu(t)+\sin\theta\mu(t))\\
&=r\left(\cos t-\cos\theta\cos(t+\theta),\sin t-\cos\theta\sin(t+\theta)\right).
\end{align*}
\par
If we take $\theta(t)=\pi/2$, a constant $\tau$ with $\cos \tau=0$ and a smooth function $\lambda: [0,2\pi) \to \R$ by $\lambda(t)=0$, then $(\beta(t)+\dot{\lambda}(t))\cos \tau-\lambda(t) \ell(t) \sin \tau=0$ for all $t \in [0,2\pi)$. 
By Definition \ref{evolutoid-involutoid}, the involutoid $\mathcal{I}nv[\tau](\gamma,\nu)$ of $(\gamma,\nu)$ is given by $\mathcal{I}nv[\tau](\gamma,\nu)(t)=\gamma(t)+\lambda(t)\mu(t)=\left(r\cos t,r\sin t\right)$. 
On the other hand, if we take $\theta(t)=\pi/2$, a constant $\tau$ with $\sin\tau=0$ and a smooth function $\lambda: [0,2\pi) \to \R$ by $\lambda(t)=-rt+c$, where $c$ is a constant, then $(\beta(t)+\dot{\lambda}(t))\cos \tau-\lambda(t) \ell(t) \sin \tau=0$ for all $t \in [0,2\pi)$. 
By Definition \ref{evolutoid-involutoid}, the involutoid $\mathcal{I}nv[\tau](\gamma,\nu)$ of $(\gamma,\nu)$ is given by 
\begin{align*}
\mathcal{I}nv[\tau](\gamma,\nu)(t)&=\gamma(t)+\lambda(t)\mu(t)\\
&=\left(r(\cos t+t\sin t)-c\sin t,r(\sin t-t\cos t)+c\cos t\right).
\end{align*}
Moreover, if we take $\theta(t)=\pi/2$ and a constant $\tau$ with $\cos\tau \sin\tau \neq0$, a smooth function $\lambda: [0,2\pi) \to \R$ by $\lambda(t)=r\cos\tau/\sin\tau+ce^{(\tan \tau)t}$, where $c$ is a constant, then $(\beta(t)+\dot{\lambda}(t))\cos \tau-\lambda(t) \ell(t) \sin \tau=0$ for all $t \in [0,2\pi)$. 
By Definition \ref{evolutoid-involutoid}, the involutoid $\mathcal{I}nv[\tau](\gamma,\nu)$ of $(\gamma,\nu)$ is given by 
\begin{align*}
\mathcal{I}nv[\tau](\gamma,\nu)(t)&=\frac{1}{\sin\tau}\bigl(r \sin\tau \cos t-r \cos\tau \sin t-c \sin\tau e^{(\tan \tau)t}\sin t,\\
&\qquad\qquad r \sin\tau \sin t+r \cos\tau \cos t+c \sin\tau e^{(\tan \tau)t} \cos t \bigr).
\end{align*}
\par
By Definition \ref{nvolute} $(1)$, $N[\theta](\gamma,\nu)$ is given by 
$$
N[\theta](\gamma,\nu)(t)=\left(r\cos t+\lambda(t)\cos(\theta+t), r\sin t+\lambda(t)\sin(\theta+t)\right),
$$
where $r+\dot{\lambda}(t)\sin \theta+\lambda(t)\cos \theta=0$ for all $t \in [0,2\pi)$. 
If $\theta=0$, then we have $\lambda(t)=-r$, $\cos(\theta+t)=\cos t$ and $\sin(\theta+t)=\sin t$. 
Therefore, we have $N[0](\gamma,\nu)(t)=\left(0, 0\right)$. 
On the other hand, if $\theta=\pi/2$, then we have $\lambda(t)=-rt+c$, $\cos(\theta+t)=-\sin t$ and $\sin(\theta+t)=\cos t$. 
Therefore, we have 
$$
N[\pi/2](\gamma,\nu)(t)=\left(r(\cos t+t\sin t)-c\sin t,r(\sin t-t\cos t)+c\cos t\right).
$$
Moreover, if we take a constant $\theta$ with $\sin\theta\cos\theta\neq0$, then we have $\lambda(t)=-r/\cos\theta+ce^{(-\cos\theta/\sin\theta)t}$.
Therefore, we have 
\begin{align*}
N[\theta](\gamma,\nu)(t)=\biggl(& r\cos t-\frac{r\cos(\theta+t)}{\cos\theta} +ce^{-\frac{\cos\theta}{\sin\theta}t}\cos(\theta+t),\\
&r\sin t-\frac{r\sin(\theta+t)}{\cos\theta}+ce^{-\frac{\cos\theta}{\sin\theta}t}\sin(\theta+t)\biggr).
\end{align*}
\par
By Definition \ref{nvolute} $(2)$, $T[\tau](\gamma,\nu)$ is given by 
$$
T[\tau](\gamma,\nu)(t)=\left(r\cos t-\lambda(t)\sin(\tau+t), r\sin t+\lambda(t)\cos(\tau+t)\right), 
$$
where $r+\dot{\lambda}(t) \cos \tau-\lambda(t) \sin \tau=0$ for all $t \in [0,2\pi)$. 
If $\tau=0$, then we have $\lambda(t)=-rt+c$, $\cos(\tau+t)=\cos t$ and $\sin(\tau+t)=\sin t$. 
Therefore, we have 
$$
T[0](\gamma,\nu)(t)=\left(r(\cos t+t\sin t)-c\sin t,r(\sin t-t\cos t)+c\cos t\right).
$$
On the other hand, if $\tau=-\pi/2$, then we have $\lambda(t)=-r$, $\cos(\tau+t)=\sin t$ and $\sin(\tau+t)=-\cos t$. 
Therefore, we have $T[-\pi/2](\gamma,\nu)(t)=\left(0, 0\right)$. 
It follows that $N[0](\gamma,\nu)$ and $T[-\pi/2](\gamma,\nu)$ are evolutes of $(\gamma,\nu)$, and $N[\pi/2](\gamma,\nu)$ and $T[0](\gamma,\nu)$ are involutes of $(\gamma,\nu)$.
Moreover, if we take a constant $\tau$ with $\sin\tau\cos\tau\neq0$, then we have $\lambda(t)=r/\sin\tau+ce^{(\sin\tau/\cos\tau)t}$.
Therefore, we have 
\begin{align*}
T[\tau](\gamma,\nu)(t)=\biggl(& r\cos t-\frac{r\sin(\tau+t)}{\sin\tau} -ce^{\frac{\sin\tau}{\cos\tau}t}\sin(\tau+t),\\
&r\sin t+\frac{r\cos(\tau+t)}{\sin\tau}+ce^{\frac{\sin\tau}{\cos\tau}t}\cos(\tau+t)\biggr).
\end{align*}
}
\end{example}
%%%%%

%%%%%
\begin{example}{\rm
Let $(\gamma,\nu):[0,2\pi) \to \R^2 \times S^1$ be 
\begin{align*}
\gamma(t)=\left(\cos^3 t, \sin^3 t\right),\ \nu(t)=\left(\sin t,\cos t\right).
\end{align*}
Then $\gamma$ is an astroid. 
By a direct calculation, $\mu(t)=(-\cos t,\sin t)$ and $(\gamma,\nu)$ is a Legendre curve with the curvature $(\ell(t),\beta(t))=(-1, 3\cos t\sin t)$. 
\par
If we take  $\theta(t)=0$, $\tau(t)=\pi/2$, a smooth function $\lambda:[0,2\pi) \to \R$ by $\lambda(t)=3\cos t\sin t$, then $\beta(t)+\lambda(t)\ell(t)=0$ for all $t \in [0,2\pi)$. 
Therefore, the Legendre curve $(\overline{\gamma},\overline{\nu})$ such that $\overline{\gamma}$ is an evolute of $(\gamma,\nu)$ is given by
\begin{align*}
\overline{\gamma}(t)&=\gamma(t)+\lambda(t)\bv(t)=\left(\cos^3 t+3\cos t\sin^2t, \sin^3t+3\cos^2t\sin t\right), \\
\overline{\nu}(t)&=-\mu(t)=\left(\cos t,-\sin t\right).
\end{align*}
\par
If we take $\theta(t)=\pi/2$, $\tau(t)=0$ and a smooth function $\lambda: [0,2\pi) \to \R$ by $\lambda(t)=(3\cos 2t)/4+c$, where $c$ is a constant, then $\beta(t)+\dot{\lambda}(t)=0$ for all $t \in [0,2\pi)$. 
Therefore, the Legendre curve $(\overline{\gamma},\overline{\nu})$ such that $\overline{\gamma}$ is an involute of $(\gamma,\nu)$ is given by
\begin{align*}
\overline{\gamma}(t)&=\gamma(t)+\lambda(t)\bv(t)=\left(-\frac{\cos^3t}{2}+\frac{3\cos t}{4}-c\cos t, -\frac{\sin^3t}{2}+\frac{3\sin t}{4}+c\sin t\right), \\
\overline{\nu}(t)&=\mu(t)=\left(-\cos t,\sin t\right).
\end{align*}
\par
If we take a constant $\theta$, $\tau(t)=\pi/2$ and a smooth function $\lambda: [0,2\pi) \to \R$, $\lambda(t)=3\cos t\sin t\cos \theta$,
then $\beta(t) \cos \theta+\lambda(t) \ell(t)=0$ for all $t \in [0,2\pi)$. 
By Definition \ref{evolutoid-involutoid}, the evolutoid $\mathcal{E}v[\theta](\gamma,\nu)$ of $(\gamma,\nu)$ is given by 
\begin{align*}
\mathcal{E}v[\theta](\gamma,\nu)(t)&=\gamma(t)+\lambda(t)(\cos\theta \nu(t)+\sin\theta \mu(t))\\
&=\left(\cos^3t+3\cos t\sin t\cos \theta \sin(t-\theta), \sin^3t+3\cos t\sin t\cos \theta \cos(t-\theta)\right).
\end{align*}
\par
If we take $\theta(t)=\pi/2$, a constant $\tau$ with $\cos \tau=0$ and a smooth function $\lambda: [0,2\pi) \to \R$ by $\lambda(t)=0$, then $(\beta(t)+\dot{\lambda}(t))\cos \tau-\lambda(t) \ell(t) \sin \tau=0$ for all $t \in [0,2\pi)$. 
By Definition \ref{evolutoid-involutoid}, the involutoid $\mathcal{I}nv[\tau](\gamma,\nu)$ of $(\gamma,\nu)$ is given by $\mathcal{I}nv[\tau](\gamma,\nu)(t)=\gamma(t)+\lambda(t)\mu(t)=\left(\cos^3 t, \sin^3 t\right)$. 
%On the other hand, if we take $\theta(t)=\pi/2$, a constant $\tau$ with $\sin\tau=0$ and a smooth function $\lambda: [0,2\pi) \to \R$ by $\lambda(t)=(3\cos2t)/4+c$, where $c$ is a constant, then $(\beta(t)+\dot{\lambda}(t))\cos \tau-\lambda(t) \ell(t) \sin \tau=0$  for all $t \in [0,2\pi)$. By Definition \ref{evolutoid-involutoid}, the involutoid $\mathcal{I}nv[\tau](\gamma,\nu)$ of $(\gamma,\nu)$ is given by \begin{align*} \mathcal{I}nv[\tau](\gamma,\nu)(t)&=\gamma(t)+\lambda(t)\mu(t)\\ &=\left(-\frac{\cos^3t}{2}+\frac{3\cos t}{4}-c\cos t, -\frac{\sin^3t}{2}+\frac{3\sin t}{4}+c\sin t\right).\end{align*}
Moreover, if we take $\theta(t)=\pi/2$, a constant $\tau$ with $\cos\tau \neq0$ and a smooth function $\lambda,: [0,2\pi) \to \R$ by
$$
\lambda(t)=-\frac{3\cos\tau(\sin2t\sin\tau-2\cos2t\cos\tau)}{2(\sin^2\tau+4\cos^2\tau)}+c e^{-(\tan \tau)t},
$$
where $c$ is a constant, then $(\beta(t)+\dot{\lambda}(t))\cos \tau-\lambda(t) \ell(t) \sin \tau=0$  for all $t \in [0,2\pi)$. 
By Definition \ref{evolutoid-involutoid}, the involutoid $\mathcal{I}nv[\tau](\gamma,\nu)$ of $(\gamma,\nu)$ is given by 
\begin{align*}
\mathcal{I}nv[\tau](\gamma,\nu)(t)&=\biggl(\cos^3t+\frac{3\cos t\cos\tau(\sin2t\sin\tau-2\cos2t\cos\tau)}{2(\sin^2\tau+4\cos^2\tau)}-c e^{-(\tan \tau)t}\cos t,\\
&\qquad\sin^3t-\frac{3\sin t\cos\tau(\sin2t\sin\tau-2\cos2t\cos\tau)}{2(\sin^2\tau+4\cos^2\tau)}+c e^{-(\tan \tau)t}\sin t\biggr).
\end{align*}
\par
By Definition \ref{nvolute} $(1)$, $N[\theta](\gamma,\nu)$ is given by
$$
N[\theta](\gamma,\nu)(t)=\left(\cos^3 t+\lambda(t)\sin(t-\theta), \sin^3 t+\lambda(t)\cos(t-\theta)\right),
$$
where $3\cos t\sin t+\dot{\lambda}(t)\sin \theta-\lambda(t)\cos \theta=0$ for all $t \in [0,2\pi)$. 
If $\theta=0$, then we have $\lambda(t)=3\cos t\sin t$, $\cos(t-\theta)=\cos t$ and $\sin(t-\theta)=\sin t$. 
Therefore, we have $N[0](\gamma,\nu)(t)=\left(\cos^3 t+3\cos t\sin^2t, \sin^3t+3\cos^2t\sin t\right)$. 
Moreover, if $\theta=\pi/2$, then we have $\lambda(t)=(3\cos 2t)/4+c$, $\cos(t-\theta)=\sin t$ and $\sin(t-\theta)=-\cos t$. 
Therefore, we have 
$$
N[\pi/2](\gamma,\nu)(t)=\left(-\frac{\cos^3t}{2}+\frac{3\cos t}{4}-c\cos t, -\frac{\sin^3t}{2}+\frac{3\sin t}{4}+c\sin t\right).
$$
Moreover, if we take a constant $\theta$ with $\sin\theta\neq0$, then we have 
$$
\lambda(t)=\frac{3(\cos\theta\sin 2t+2\sin\theta\cos2t)}{2(4\sin^2\theta+\cos^2\theta)}+ce^{\frac{\cos\theta}{\sin\theta}t}.
$$
Therefore, we have 
\begin{align*}
N[\theta](\gamma,\nu)(t)=\biggl(& \cos^3t+\frac{3(\cos\theta\sin 2t+2\sin\theta\cos2t)}{2(4\sin^2\theta+\cos^2\theta)}\sin(t-\theta)+ce^{\frac{\cos\theta}{\sin\theta}t}\sin(t-\theta),\\
&\sin^3t+\frac{3(\cos\theta\sin 2t+2\sin\theta\cos2t)}{2(4\sin^2\theta+\cos^2\theta)}\cos(t-\theta)+ce^{\frac{\cos\theta}{\sin\theta}t}\cos(t-\theta)\biggr).
\end{align*}
By Definition \ref{nvolute} $(2)$, $T[\tau](\gamma,\nu)$ is given by 
$$
T[\tau](\gamma,\nu)(t)=\left(\cos^3 t-\lambda(t)\cos(t-\tau), \sin^3 t+\lambda(t)\sin(t-\tau)\right), 
$$
where $3\cos t\sin t+\dot{\lambda}(t) \cos \tau+\lambda(t) \sin \tau=0$ for all $t \in [0,2\pi)$. 
If $\tau=0$, then we have $\lambda(t)=(3\cos 2t)/4+c$, $\cos(t-\tau)=\cos t$ and $\sin(t-\tau)=\sin t$. 
Therefore, we have 
$$
T[0](\gamma,\nu)(t)=\left(-\frac{\cos^3t}{2}+\frac{3\cos t}{4}-c\cos t, -\frac{\sin^3t}{2}+\frac{3\sin t}{4}+c\sin t\right).
$$
Moreover, if $\tau=-\pi/2$, then we have $\lambda(t)=3\cos t\sin t$, $\cos(t-\tau)=-\sin t$ and $\sin(t-\tau)=\cos t$. 
Therefore, we have $T[-\pi/2](\gamma,\nu)(t)=\left(\cos^3 t+3\cos t\sin^2t, \sin^3t+3\cos^2t\sin t\right)$. 
It follows that $N[0](\gamma,\nu)$ and $T[-\pi/2](\gamma,\nu)$ are evolutes of $(\gamma,\nu)$, and $N[\pi/2](\gamma,\nu)$ and $T[0](\gamma,\nu)$ are involutes of $(\gamma,\nu)$.
Moreover, if we take a constant $\tau$ with $\cos\tau\neq0$, then we have 
$$
\lambda(t)=-\frac{3(\sin2t\sin\tau-2\cos2t\cos\tau)}{2(\sin^2\tau+4\cos^2\tau)}+c e^{-(\tan \tau)t}.
$$
Therefore, we have 
\begin{align*}
T[\tau](\gamma,\nu)(t)=\biggl(& \cos^3t+\frac{3(\sin2t\sin\tau-2\cos2t\cos\tau)}{2(\sin^2\tau+4\cos^2\tau)}\cos(t-\tau)-ce^{-\tan \tau}\cos(t-\tau),\\
&\sin^3t-\frac{3(\sin2t\sin\tau-2\cos2t\cos\tau)}{2(\sin^2\tau+4\cos^2\tau)}\sin(t-\tau)+ce^{-\tan \tau}\sin(t-\tau)\biggr).
\end{align*}
}
\end{example}
%%%%%

%%%%%%%%%% Bibliography %%%%%%%%%%%%%%%%%%

%%%%%%%%
Nozomi Nakatsuyama, 
\\
Muroran Institute of Technology, Muroran 050-8585, Japan,
\\
E-mail address: 25096009b@muroran-it.ac.jp
\\
\\
Masatomo Takahashi, 
\\
Muroran Institute of Technology, Muroran 050-8585, Japan,
\\
E-mail address: masatomo@muroran-it.ac.jp


\begin{thebibliography}{99}
{\small

\bibitem{Aminov} Y. Aminov, 
\newblock{Differential geometry and the topology of curves. Translated from the Russian by V. Gorkavy.} 
\newblock{Gordon and Breach Science Publishers}, Amsterdam, (2000).

\bibitem{AM} T. M. Apostol, M. A. Mnatsakanian,
\newblock{Tanvolutes: generalized involutes}.
\newblock{Amer. Math. Monthly} {\bf 117}, (2010), 701--713.

%\bibitem{Arnold1} V. I. Arnol'd,\newblock{\em Singularities of Caustics and Wave Fronts}.\newblock Mathematics and Its Applications {\bf 62} Kluwer Academic Publishers, (1990).

%\bibitem{Arnold2} V. I. Arnol'd, S. M. Gusein-Zade and A. N. Varchenko,\newblock{\em Singularities of Differentiable Maps vol. {\rm I}}. \newblock Birkh\"auser, (1986). 

\bibitem{Banchoff-Lovett} Y. Banchoff, S. Lovett, 
\newblock{Differential geometry of curves and surfaces.}
\newblock{A K Peters, Ltd., Natick,} MA, (2010). 

\bibitem{Berger-Gostiaux} M. Berger, B. Gostiaux,  
\newblock{Differential geometry: manifolds, curves, and surfaces. Translated from the French by Silvio Levy.} 
\newblock{Graduate Texts in Mathematics, 115. Springer-Verlag,} New York, (1988).

\bibitem{Bertrand} J. Bertrand, 
\newblock{M\'emoire sur la th\'eorie des courbes \`a double courbure.}
\newblock{J. de meth\'ematiques pures et appliqu\'ees}. {\bf 15}, (1850), 332--350.

%\bibitem{Bishop} R. L. Bishop, \newblock{There is more than one way to frame a curve.}\newblock{American Mathematical Monthly}. {\bf 82} (1975), 246--251.

\bibitem{Bruce-Gaffney} J. W. Bruce, T. J. Gaffney, 
\newblock{Simple singularities of mappings {${C},0\rightarrow {C}^{2},0$}}, 
\newblock{J. London Math. Soc. (2).} {\bf 26}, (1982), 465--474.

%\bibitem{Bruce-Giblin} J. W. Bruce, P. J. Giblin, \newblock{Curves and singularities. A geometrical introduction to singularity theory. Second edition}. \newblock{Cambridge University Press}, Cambridge, (1992).

\bibitem{doCarmo} M. P. do Carmo, 
\newblock{Differential geometry of curves and surfaces.} 
Translated from the Portuguese. Prentice-Hall, Inc., Englewood Cliffs, N.J., (1976).

\bibitem{Fukunaga-Takahashi-2013} T. Fukunaga and M. Takahashi,
\newblock{Existence and uniqueness for Legendre curves}. 
\newblock{J. Geom.} {\bf 104}, (2013), 297--307.

\bibitem{Fukunaga-Takahashi-2015} T. Fukunaga, M. Takahashi,
\newblock{Evolutes and involutes of frontals in the Euclidean plane}.
\newblock{Demonstr. Math.} {\bf 48}, (2015), 147--166.

%\bibitem{Fukunaga-Takahashi-2017} T. Fukunaga, M. Takahashi, \newblock{Existence conditions of framed curves for smooth curves.} \newblock{Journal of Geometry.} {\bf 108}, (2017), 763--774. doi: 10.1007/s00022-017-0371-5 

%\bibitem{Honda-Takahashi-2016} S. Honda, M. Takahashi, \newblock{Framed curves in the Euclidean space.}\newblock{Advances in Geometry.} {\bf 16}, (2017), 265--276.doi: 10.1515/advgeom-2015-0035

\bibitem{Honda-Takahashi-2020} S. Honda, M. Takahashi, 
\newblock{Bertrand and Mannheim curves of framed curves in the 3-dimensional Euclidean space.}
\newblock{Turkish J. Math.} {\bf 44}, (2020), 883--899.

%\bibitem{Honda-Takahashi-Preprint} S. Honda, M. Takahashi, \newblock{Circular evolutes and involutes of framed curves in the 3-dimensional Euclidean space.}\newblock{Preprint}, (2024).

\bibitem{Giblin-Warder} P. J. Giblin, J. P. Warder, 
\newblock{Evolving evolutoids}. 
\newblock{Amer. Math. Monthly} {\bf 121}, (2014), 871--889.

\bibitem{Gibson} C. G. Gibson,
\newblock{Elementary geometry of differentiable curves. An undergraduate introduction}.
\newblock{Cambridge University Press}, Cambridge, (2001).

\bibitem{Gray} A. Gray, E. Abbena and S. Salamon, 
\newblock{Modern differential geometry of curves and surfaces with Mathematica. Third edition.} 
\newblock{Studies in Advanced Mathematics. Chapman and Hall/CRC}, Boca Raton, FL, (2006).

\bibitem{HCIP} J. Huang, L. Chen, S. Izumiya, D. Pei, 
\newblock{Geometry of special curves and surfaces in 3-space form.}
\newblock{J. Geom. Phys.} {\bf 136}, (2019), 31--38.

%\bibitem{Ishikawa-book} G. Ishikawa, \newblock{\em Singularities of Curves and Surfaces in Various Geometric Problems}.\newblock{CAS Lecture Notes 10, Exact Sciences}. (2015).

%\bibitem{Izumiya09} S, Izumiya.  Legendrian dualities and spacelike hypersurfaces in the lightcone.Moscow Mathmatical Journal. {\bf 9}, (2009), 325-357. https://doi.org/10.17323/1609-4514-2009-9-2-325-357

%\bibitem{Izumiya-book} S. Izumiya, M. C. Romero-Fuster, M. A. S. Ruas and F. Tari,\newblock{\em Differential Geometry from a Singularity Theory Viewpoint}.\newblock{World Scientific Pub. Co Inc.} (2015).

\bibitem{Izumiya-Takeuchi1} S. Izumiya, N. Takeuchi,  
\newblock{Generic properties of helices and Bertrand curves.}
\newblock{J. Geom.} {\bf 74}, (2002), 97--109.

\bibitem{Izumiya-Takeuchi-2019} S. Izumiya, N. Takeuchi, 
\newblock{Evolutoids and pedaloids of plane curves}.
\newblock{Note Mat.} {\bf 39}, (2019), 13--23.

\bibitem{Kuhnel} W. K\"uhnel, 
\newblock{Differential geometry. Curves-surfaces-manifolds. Translated from the 1999 German original by Bruce Hunt.} 
\newblock{Student Mathematical Library, 16. American Mathematical Society,} Providence, RI, (2002).

\bibitem{Li-Pei} E. Li, D. Pei, 
\newblock{Enveloids and involutoids of spherical Legendre curves}. 
\newblock{J. Geom. Phys.} {\bf 170}, (2021), Paper No. 104371, 23 pp.

\bibitem{Liu-Wang} H. Liu, F. Wang, 
\newblock{Mannheim partner curves in 3-space.}
\newblock{J. Geom.} {\bf 88}, (2008), 120--126.

%\bibitem{Lucas-Ortega} P. Lucas, J. A. Ortega-Yag\"ues, \newblock{A variational characterization and geometric integration for Bertrand curves.} \newblock{Journal of Mathematical Physics} {\bf 54}, 2013, 043508, 12 pp.

\bibitem{Nakatsuyama} N. Nakatsuyama,
\newblock{Evolutes and involutes of framed curves in the Euclidean 3-space.}
\newblock{In preparation}, (2026).

%\bibitem{Nakatsuyama-Takahashi1} N. Nakatsuyama, M. Takahashi, \newblock{On vertices of frontals in the Euclidean plane.}  \newblock{Bull. Braz. Math. Soc. (N.S.)} {\bf 55}, (2024), Paper No. 35, 21 pp.

\bibitem{Nakatsuyama-Takahashi2} N. Nakatsuyama, M. Takahashi, 
\newblock{Bertrand types of regular curves and Bertrand framed curves in the Euclidean 3-space} 
\newblock{To appear in Hokkaido Mathematical Journal}, (2026).

%\bibitem{Nakatsuyama-Takahashi} N. Nakatsuyama, M. Takahashi, \newblock{Bertrand lightcone framed curves in the Lorentz-Minkowski 3-space.} \newblock{In preparation}, (2024).

\bibitem{Papaioannou-Kiritsis} S. G. Papaioannou, D. Kiritsis,
\newblock{An application of Bertrand curves and surfaces to CADCAM.} 
\newblock{Computer-Aided Design}. {\bf 17}, (1985), 348--352.

\bibitem{Porteous} I. R. Porteous,  
\newblock{Geometric differentiation for the intelligence of curves and surfaces, Second edition.} \newblock{Cambridge University Press,} Cambridge, (2001).

\bibitem{Struik} D. J. Struik, 
\newblock{Lectures on classical differential geometry. Reprint of the second edition.} 
\newblock{Dover Publications, Inc.,} New York, (1988).

%\bibitem{Takahashi} M. Takahashi, \newblock{Legendre curves in the unit spherical bundle over the unit sphere and evolutes.} \newblock{Contemp. Math.} {\bf 675}, (2016), 337--355.

}
\end{thebibliography}
\end{document}